\documentclass[journal,twocolumn]{IEEEtran}

\ifCLASSINFOpdf
\usepackage[pdftex]{graphicx}
\else
\usepackage[dvips]{graphicx}
\fi

\usepackage{amsmath, amssymb, amsfonts}
\usepackage{bbm} 
\usepackage{algorithmicx}
\usepackage{algorithm}
\usepackage{algpseudocode}
\ifCLASSOPTIONcompsoc
  \usepackage[caption=false,font=normalsize,labelfont=sf,textfont=sf]{subfig}
\else
  \usepackage[caption=false,font=footnotesize]{subfig}
\fi
	
	\usepackage{longtable}

\begin{document}

		\title{
			Evaluation of Prosumer Networks for Peak Load Management in Iran: A Distributed Contextual Stochastic Optimization Approach  
		}
		
		\author{Amir Noori, Babak Tavassoli, and Alireza Fereidunian
			\thanks{Amir Noori, Babak Tavassoli, and Alireza Fereidunian are with the Faculty of Electrical Engineering, KN Toosi of Technology, Tehran,
				Iran, e-mail:} 
		}
		
	\markboth{Journal of \LaTeX\ Class Files,~Vol.~14, No.~8, Jun~2024}%
	{Shell \MakeLowercase{\textit{et al.}}: Evaluation of Prosumer Networks for Peak Load Management in Iran: A Distributed Contextual Stochastic Optimization Approach}
	\maketitle
	
	\begin{abstract}
		Renewable prosumers face the complex challenge of balancing self-sufficiency with seamless grid and market integration. This paper introduces a novel prosumers network framework aimed at mitigating peak loads in Iran, particularly under the uncertainties inherent in renewable energy generation and demand.
		A cost-oriented integrated prediction and optimization approach is proposed, empowering prosumers to make informed decisions within a distributed contextual stochastic optimization (DCSO) framework. 
		The problem is formulated as a bi-level two-stage multi-time scale optimization to determine optimal operation and interaction strategies under various scenarios, considering flexible resources.
		To facilitate grid integration, a novel consensus-based contextual information sharing mechanism is proposed. This approach enables coordinated collective behaviors and leverages contextual data more effectively. 
		The overall problem is recast as a mixed-integer linear program (MILP) by incorporating optimality conditions and linearizing complementarity constraints. Additionally, a distributed algorithm using the consensus alternating direction method of multipliers (ADMM) is presented for computational tractability and privacy preservation. Numerical results highlights that integrating prediction with optimization and implementing a contextual information-sharing network among prosumers significantly reduces peak loads as well as total costs.
	\end{abstract}
	

	\IEEEpeerreviewmaketitle
	
	\begin{table}[htbp]
		\centering
		\label{tab:notation}
		\begin{tabular}{p{0.5cm}p{7.5cm}}
			\multicolumn{2}{l}{\textbf{Notation \quad Description}} \\
			\multicolumn{2}{l}{\textit{Sets and Indices}} \\
			$m, n$ & Indices for prosumers \\
			$i, j$ & Indices for data samples \\
			$k$ & Number of neighbor points in kNN \\
			$h$ & Index for time periods \\
			$\nu$ & Index for iterations \\
			$\mathcal{N}, \mathcal{N}_n$ & Set of prosumers, neighbors of prosumers $n$ \\
			$\mathcal{H}$ & Set of time horizons \\
			\multicolumn{2}{l}{Parameters} \\
			${p}_n^g$ & Photovoltaic (PV) power output of prosumer $n$ (kW) \\
			${p}_n^l$ & Power of must-run loads of prosumer $n$ (kW) \\
			$p_n^{s,r}$ & Desired power consumption of shiftable loads (SLs) (kW) \\
			$C_n^s$ & Total SL energy demand within horizon $\mathcal{H}$ (kW) \\
			$\alpha, \beta$ & Coefficients of cost functions \\
			$\eta_n$ & Efficiency of Battery Energy Storage System (BESS) \\
			$c_{nm}$ & Energy price between prosumer $n$ and $m$ (\textcent/kWh) \\
			$\Delta t$ & Time step size (hour) \\
			$S$ & Number of samples \\
			\multicolumn{2}{l}{\textit{Decision Variables}} \\
			$p_n^c, p_n^c$ & Charging/discharging power of BESS (kW) \\
			$p_n^b$ & Power of BESS (kW) \\
			$p_n^s$ & Consumption power of SLs (kW) \\
			$p_{n}^{mt}$ & Day-ahead Energy trade of prosumer $n$ with the main grid (kWh) \\
			$q_{n}^{mt}$ & Real-time Energy trade of prosumer $n$ with the main grid (kWh) \\
			$p_{nm}$ & P2P Energy trade of prosumer $n$ and $m$ (kWh) \\
			$p_n$ & Total P2P energy trade of prosumer $n$ (kWh) \\
			$e_{n}$ & State-of-Charge (SoC) of BESS of prosumer $n$ (kWh) \\
			$s_{n}$ & Energy shift state of SLs of prosumer $n$ (kWh) \\
		\end{tabular}
	\end{table}
	
	\section{Introduction}\label{sec1}

	\subsection{Motivation and Background} \label{sec1.1}
	\IEEEPARstart{I}{ran} experiences significant electricity peak loads during summer, when cooling demands cause consumption frequently exceeds production capacity. This seasonal surge in demand places immense pressure on the national grid, often resulting in power shortages, particularly in the residential and industrial sectors. Expanding production capacity by an estimated 14,000 to 17,000 MW \cite{Qasemi2024} to meet these peaks is not economically viable, given the short duration of peak periods.
	
	As a strategic alternative, the implementation of prosumer networks offers an efficient solution. Prosumers can coordinate their generation and consumption behaviors to help balance supply and demand more effectively, thereby reducing the strain on the grid during peak periods. These networks not only address the immediate issue of peak load mitigation on a large scale but also accelerate grid modernization.
	
	The flexibility offered by energy storage systems, such as batteries, further strengthens the grid's ability to absorb and dispatch renewable energy, smoothing out fluctuations in supply and demand. Essentially, prosumer networks promote the integration of renewable energy sources and the deployment of flexible resources, contributing to a more sustainable and robust energy infrastructure.
	
	Moreover, prosumer networks foster a more decentralized and resilient energy system. By distributing generation and storage across a wide array of small-scale producers, the grid becomes less vulnerable to large-scale outages and more adaptable to changes in energy consumption patterns. This decentralized approach also empowers communities and individuals, allowing them to take an active role in energy management and contribute to the overall sustainability of the power system.
	
	Beyond these technical and environmental benefits, prosumer networks can drive economic growth by creating new opportunities for investment in renewable energy, smart grid technologies, and digitalization. The development of such networks supports the expansion of renewable energy infrastructure in Iran, promoting innovation and creating jobs in the clean energy sector.

	\subsection{Related Works}
	
	The EU's Clean Energy Package introduced provisions for individual and collective self-consumption (CSC) and renewable energy communities (RECs), empowering citizens to engage in the energy market and fostering decentralized renewable energy production \cite{capros2018outlook}. For effective adoption, targeted measures such as information dissemination, facilitation, and open participation are crucial. A comparative study highlights the variations in these frameworks across European countries, examining factors like geographical limitations, grid tariffs, and membership requirements in energy communities \cite{frieden2020collective}.
	
	Globally, similar policies, including net metering, feed-in tariffs, and tax incentives, have been implemented or are under development. However, Iran faces unique challenges in adopting these policies due to economic difficulties that have hampered the development of new power generation, particularly renewable energy, leading to significant energy imbalances during peak summer demand \cite{Qasemi2024}. While efforts to develop small-scale power plants and renewable projects in remote areas are well-intentioned, they have not adequately addressed the growing supply-demand gap. A more comprehensive strategy is needed to ensure long-term energy security in Iran.
	
	Demand-side management (DSM) has been a cornerstone of peak load mitigation for decades \cite{gellings2017evolving}. Price-based mechanisms, such as dynamic pricing \cite{yang2018real} and time-of-use tariffs, and incentive-based approaches, like demand response programs \cite{albadi2007demand}, are common DSM strategies. DSM research encompasses a diverse range of participants, energy sources, and market configurations \cite{shakeri2020overview}. Recent technological advancements and the increasing integration of renewable energy have created an enabling environment for federated energy systems \cite{morstyn2018using}, empowering various stakeholders to actively participate in grid operations.
	
	The exploration of prosumer dynamics within contemporary energy systems has advanced through various studies, each contributing to the understanding of how prosumers can support flexibility and grid stability beyond mere cost reduction. One review challenges traditional views, offering insights into aggregation models and the complexities of coordinating prosumer services with grid operators \cite{grvzanic2022prosumers}. Another study evaluates the economic implications of complementarity, cost allocation, and demand-response programs, demonstrating their potential to lower operational costs in energy communities through mixed-integer linear programming \cite{volpato2022general}.
	
	A decentralized framework that coordinates prosumer flexibility scheduling, aimed at mitigating peak loads, has been proposed, incorporating a multi-objective optimization approach and collective learning algorithms \cite{mashlakov2021decentralized}. Comparative analysis of regulatory frameworks across multiple EU countries reveals conditions favorable to collective renewable energy prosumers, with recommendations for policy improvements \cite{ines2020regulatory}. The concept of the "co-prosumer," emphasizing collective action and resilience, is applied to energy and agriculture, proposing a model that balances individual and collective interests \cite{ritzel2022prosuming}.
	
	Aggregation policies such as virtual metering and peer-to-peer (P2P) trading are also explored, highlighting their role in enhancing distributed generation and market access \cite{moura2019prosumer}. Additionally, optimization models like mixed-integer second-order conic programming (MISOCP) are developed for load and energy management in active distribution networks, integrating renewable technologies and demand-side management \cite{alamolhoda2024integrated}. Another bi-level program addresses the interaction between prosumers and the grid, focusing on maximizing self-consumption and minimizing generation costs \cite{riaz2017generic}.

	The literature addressing uncertainties in dynamic environments between large prosumers spans a wide range of models, from stochastic and robust optimization to machine learning techniques. A hybrid stochastic and robust optimization approach, as proposed in \cite{fanzeres2014contracting}, optimizes contracting strategies for renewable energy trading companies, considering market uncertainties. Similarly, a binary prediction market designed in \cite{shamsi2021prediction} enhances probabilistic renewable energy forecasts, thereby reducing electricity market imbalance costs and improving decision-making under uncertainty.
	Stochastic optimization models also play a significant role, as demonstrated by \cite{do2021stochastic}, who present a model for day-ahead market clearing and settlement. This model employs Monte Carlo simulations to account for distributed energy resource (DER) uncertainties, effectively minimizing expected system operation costs. On the data-driven front, \cite{stratigakos2022prescriptive} integrate forecasting and optimization into a single module using decision trees and sample average approximation, improving trading decisions in renewable energy markets over standard stochastic approaches.
	By using principal component analysis-based approximation to manage uncertainties, their model delivers cost savings and increased renewable penetration. Distributed algorithms are also explored, with \cite{zhou2019online} developing a distributed stochastic dual gradient algorithm for DER management in real-time settings, ensuring guaranteed convergence and efficient coordination. Additionally, \cite{liu2022distributed} introduce a primal-dual-based distributed algorithm for managing multi-objective energy systems, balancing costs, emissions, and flexibility while reducing communication burdens.
	
	Recently, contextual optimization methods have emerged as a powerful tool for decision-making under uncertainty. A survey by \cite{sadana2024survey} explores these methods, which integrate prediction algorithms with optimization techniques. A bi-level framework proposed by \cite{munoz2022bilevel} further exemplifies this approach, optimizing decision-making with contextual information to maximize decision value while ensuring feasibility. In addition, \cite{stratigakos2022end} introduce an end-to-end learning framework using prescriptive trees, which enhances hierarchical renewable energy production forecasts and mitigates the impact of missing data, outperforming traditional two-step reconciliation methods.
	
	This paper introduces a distributed contextual stochastic optimization (DCSO) framework to enable prosumers to collaborate effectively in energy management, incorporating P2P trading, flexible resource coordination, and a contextual information marketplace.
		
	\subsection{Contributions}
	The contributions of this research article can be summarized as follows:
	
	\begin{itemize}
		\item[--] Contextual Stochastic Framework: Developed a framework integrating prediction and optimization for prosumers in energy markets.
		
		\item[--] Distributed Analytics: Introduced a distributed context-sharing mechanism and prescription algorithm for enhanced coordination.
		
		\item[--] Practical Implementation: Utilized distributed algorithms to ensure privacy and scalability, enabling real-world deployment.
	\end{itemize}
	
	The rest of this paper is organized as follows. Section \ref{sec2}, models the energy system of the proposed energy community. This is followed by a description of the contextual stochastic optimization framework to model uncertainties and covariates in Section \ref{sec3}. Simulation results are presented in Section \ref{sec4} and Section \ref{sec5} concludes the paper.

	\textit{Notation}. Operator $[a]^+ := \max\{a, 0\}$ is the projection to the non-negative reals. The expectation is denoted by $\mathbb{E}[\cdot]$. The indicator function of set $A$ is denoted by $\mathbbm{1}_A(x)$, which equals $1$ if $x\in A$ and $0$ otherwise.

	\section{Problem Formulation}\label{sec2}
	
	This section presents the system model and prosumer model that form the foundation of our research, providing a framework for analyzing the interactions between prosumers and the grid.

	\subsection{System Model}
	
	Renewable-based energy communities encounter significant challenges due to the inherent uncertainties in energy generation and consumption. To address these challenges, P2P energy trading has been proposed as a method to enhance energy efficiency and local consumption. However, the reliance on renewable resources for P2P trading without coordination with other prosumers and network can exacerbate imbalances within the distribution system.
	
\begin{figure}[!t]
		\centering
		\includegraphics[width=\columnwidth]{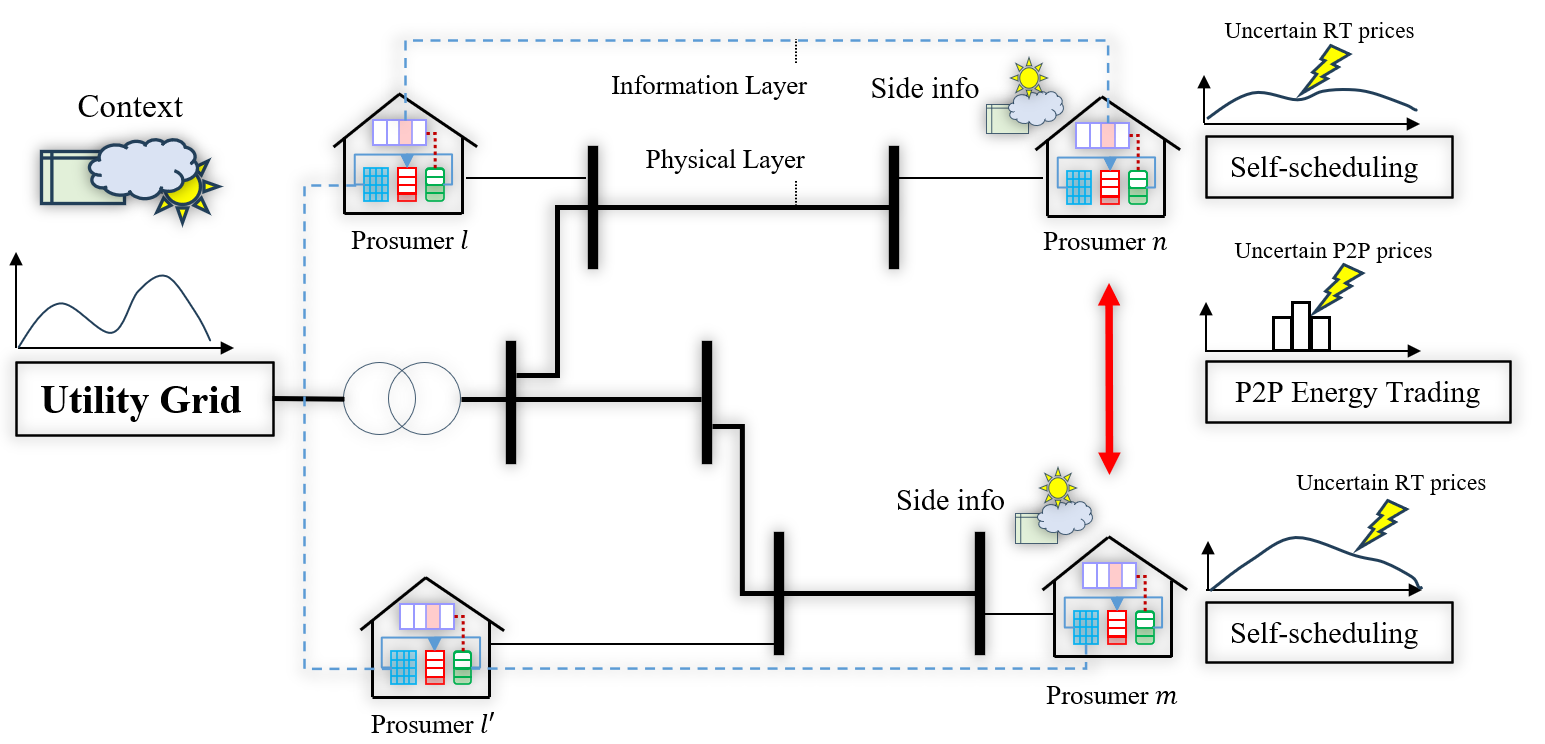}
		\caption{Conceptual overview of the proposed distributed, context-aware framework for coordinated peer-to-peer energy trading and flexibility provision.}
		\label{fig1}
	\end{figure}

	Fig. \ref{fig1} illustrates the imbalances that may arise when prosumers engage in P2P energy transactions. For instance, prosumer $n$ injects power $p_{nm}$ into bus $i$, while prosumer $m$ withdraws power $p_{mn}$ from bus $j$. Such exchanges can cause issues like line congestion, requiring intervention to manage excess power flow or voltage imbalances. While traditional methods involve the distribution system operator or other parties to handle these imbalances, this study explores a prosumer-based approach. Each prosumer $n\in\mathcal{N}$ is assumed to have bilateral trade agreements with energy partners $\mathcal{N}_n$. These agreements, managed via Home Energy Management Systems (HEMS), aim to maintain reliable and efficient power balance during P2P transactions. By combining energy trading with flexibility resources, prosumers can mitigate the impacts of uncertain generation and consumption, ensuring stable and efficient operation within the distribution system.
	
	\subsection{Prosumer Model}
	
	In this model, prosumers play a dual role as both consumers and producers. They generate their own electricity, manage their energy storage, and participate in the energy market. To ensure a balance between supply and demand, each prosumer $n$ must meet their energy consumption needs.

	\begin{equation}
		p_n^{mt} + q_n^{mt} + p_n^g + p_n = p_n^l  + p_n^b + p_n^s,\\  \label{powerbalance} 
	\end{equation}
	
	The intermittent generation of renewable power $p^g$ and fluctuating demand $p^l$ create disruptions in the power balance for prosumers, impacting the distribution network. Strategies for re-balancing include leveraging grid resources, local energy storage systems such as battery energy storage systems (BESS) and shiftable loads (SL), as well as P2P energy trading. However, prosumers exhibit preferences in adopting these strategies. The degradation cost of a BESS, a primary cost driver, can be effectively modeled using a quadratic function \cite{noori2024joint}
	\begin{align*}
		& C_n^{b} (p_{n}^{b}) = \alpha_n^{b} \lVert p_{n}^{b} \rVert ^2 
	\end{align*}
	where $\alpha^{b}$ and $\beta^{b}$ are positive weights. The operation of the BESS adheres to the constraints defined by the following equations
	\begin{subequations}\label{st}
		\begin{align}
			e_{n,h+1} &= e_{n,h} + \eta_n p_{n,h}^{b} 	\label{st1} \\  
			\underline{e}_n &\leq  e_{n,h+1} \leq \bar{e}_n 	\label{st2} \\
			0 &\leq p_{n,h}^{b} \leq \bar{p}_{n}^{b} 	\label{st3}
		\end{align}
	\end{subequations}
	
	Users with SLs have energy consumption profiles tailored to their living habits. Deviations from these preferred profiles can inconvenience end-users, quantified as discomfort costs \cite{noori2024joint}.
	\begin{align*}
		C_{n}^{s}(p_{n}^{s})= \alpha _n^{s} \lVert s_{n} \rVert ^2 
	\end{align*}
	where $\alpha^{s}$ and $\beta^{s}$ are positive weights.
	
	The operation of SLs is governed by the following equations
	\begin{subequations}\label{sh}
		\begin{align}
			s_{n,h+1} &= s_{n,h} + (p_{n}^{s} - p_{n}^{s,r}) \\
			\sum _{h \in \mathcal{H}} p_{n,h}^{s} \Delta h &= C_n^s 	\label{sh1} \\
			-\bar{p}_{n}^{s} &\leq p_{n}^{s} \leq \bar{p}_{n}^{s} \label{sh2}
		\end{align}
	\end{subequations}
	where $s$ denotes the energy shift state of SLs. 
	
	Prosumers must make day-ahead energy trading decisions before random parameters are realized to enable optimal resource allocation for both themselves and the main grid. 
	The objective of prosumers in energy trading with the main grid and other prosumers can be formulated as follows
	\begin{align*}
		C_n^{e}(p_n^{mt}, p_{nm})  = (c_p^{mt})^{\top} p_{n}^{mt} + (\tilde{c}_q^{mt})^{\top} q_{n}^{mt} + (\tilde{c}_{nm})^{\top} \, p_{nm}
	\end{align*}
	where the first and second terms represent the day-ahead and real-time purchase of energy from the main grid, and the last term represents energy purchased from other prosumers.
	Maintaining mutual energy balance between prosumers and trading excess energy are ensured by the following set of equations
	\begin{subequations}\label{recipe}
		\begin{align}
			& p_{nm} + p_{mn}=0, && \forall m \in \mathcal{N}_n \label{recipe1} \\
			& -p_{n}^g -p_{n}^{d} \leq p_{nm}^{e} \leq p_{n}^{c} + p_{n}^{s}, && \forall m \in \mathcal{N}_n \label{recipe2}
		\end{align}
	\end{subequations}
	The inequality \eqref{recipe2} ensures that only locally generated energy is traded between prosumers, thereby preventing energy arbitrage.
	To simplify the exposition, the total decision variables will henceforth be denoted by $z_n=\{p_{n}^b, p_{n}^s, p_n^{mt}, p_{nm}\}$, uncertain parameters as $Y_n=\{\tilde{c}_q^{mt}, \tilde{c}_{nm}\}$ and the feasible solution set as $\mathcal{Z}_n=\{z_n | \, \text{Eqs.} \, \eqref{powerbalance}, \eqref{st}, \eqref{sh}, \eqref{recipe} \}$.

	\begin{figure}[!t]
		\centering
		\includegraphics[width=\columnwidth]{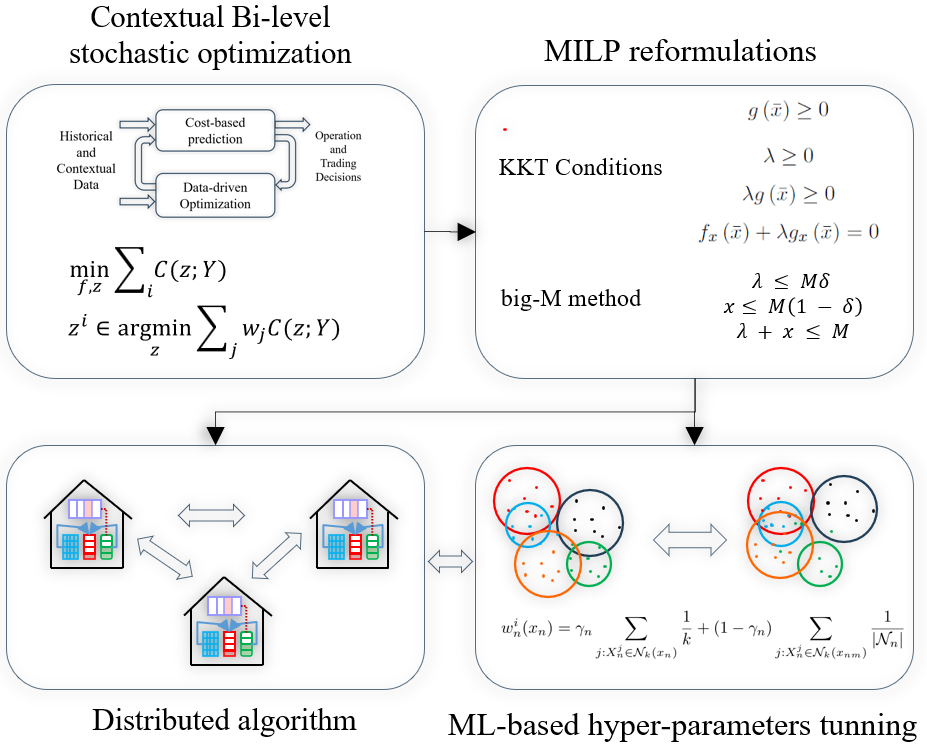}
		\caption{Distributed Contextual Stochastic Optimization Framework}
		\label{fig2}
	\end{figure}

	\section{Decision-Making under Uncertainty}\label{sec3}
	
	Prosumers face an optimization challenge characterized by uncertain parameters within a data-rich environment, necessitating a data-driven stochastic optimization approach. This section formulates the contextual stochastic optimization problem by integrating contextual data and introduces a novel distributed stochastic optimization framework.
	
	\subsection{Contextual Stochastic Optimization}
	
	Stochastic optimization in power systems is frequently challenged by uncertainties in both the objective functions and system constraints. To streamline the analysis and isolate the effects of uncertain objectives, this study focuses exclusively on cases where uncertainty is confined to the objective function. Extending the proposed method to account for uncertainty in constraints through the use of chance constraints or robust optimization techniques is a well-established approach in the literature and can be readily integrated into the framework \cite{noori2024joint}\cite{}. Prosumers within a stochastic optimization framework seek to minimize the following convex optimization problem.
	\begin{align*}
		\min_{z_n \in \mathcal{Z}_n} \mathbb{E} [C(z_{n};{Y}_n)]    
	\end{align*}	
	where the expectation is computed over the distribution of uncertain parameter $Y$. A more generalized approach integrates contextual information $X$ into the decision-making problem.
	\begin{align*}
		\min_{z_n \in \mathcal{Z}_n} \mathbb{E}[C(z_{n};Y_n) \, \arrowvert \, X_n=x_n]
	\end{align*}
	Here, the decision-maker seeks the optimal action $z_n^*(x)\in \mathcal{Z}_n$ that minimizes the expected costs conditioned on the context $x$, as expressed through the conditional distribution $Q_{y \arrowvert x}$ of $y$ given $x$
	\begin{align} \label{CSO}
		\min_{z_n \in \mathcal{Z}_n} \mathbb{E}_{Q_{y\arrowvert x}}[C(z_{n};Y_n)]
	\end{align}
	
	In the subsequent sections, the problem is reformulated into a tractable optimization model by integrating the prediction of uncertain parameters using covariate data.
	
	\subsection{Bi-level Integrated Prediction and Optimization}
	Addressing the challenges posed by problems like the above is complex due to the often unknown probability distribution. However, historical data on uncertain parameters and covariates-represented as ${(X_n^1,Y_n^1),\dots, (X_n^S,Y_n^S)}$ is usually available. The sample average approximation (SAA) method \cite{kleywegt2002sample} is fundamental in data-driven stochastic optimization, leveraging the empirical distribution derived from the samples to approximate the true expectation.
	
	Given a new covariate $x_n$, the conditional expectation \eqref{CSO} can be estimated using the following estimator as described in \cite{bertsimas2020predictive}
	\begin{align*}
		\min_{z_n \in \mathcal{Z}_n} \sum_{i=1}^S w_n^i(x_n) C(z_{n};Y_n^i)
	\end{align*}
	The weights $w_n^i(x_n)$ are non-negative and sum to one, and are critical in maintaining consistency and asymptotic optimality.
	
	\begin{align*}
		& w_n^i(x_n) \geq 0, \quad \forall \, 1 \leq i \leq n, \quad
		\sum_{i=1}^S w_n^i(x_n) = 1 
	\end{align*}
	These weights are typically derived from the non-parametric function $f: \mathbb{R}^{d_x} \rightarrow \mathbb{R}^{d_y}$ based on historical training data and the observed covariate $x_n$. 
	As detailed in \cite{bertsimas2020predictive}, an integrated framework to optimize for the function $f$ that leads to good perscriptive performance of its decisions can be given by 
		\begin{align}\label{CSOKNN}
			&\min_{f \in \mathcal{F}, z_n(f,X_n^i)} \sum_{i=1}^S C(z_n(f,X_n^i); Y_n^i) \nonumber \\
			&\text{subject to} \\
			&z_n(f,X_n^i) \in \arg\min_{z_n\in\mathcal{Z}_n} \sum_{j=1}^S w_n^j(X_n^i) C(z_n^i; Y_n^j), \forall 1 \leq i \leq S. \nonumber
		\end{align}

	\subsection{Distributed Model Parameters Tuning}	
	Various classical machine learning algorithms such as k-nearest neighbors (kNN), Nadaraya–Watson kernel regression, and Random Forest have been suggested in the literature \cite{mundru2019predictive} for computing the weights $w_n^i(x)$.
	For example, in the standard kNN method, weights are computed as
	\begin{align*}
		w_n^i(x_n) = 
		\sum_{j=1}^S \frac{1}{k} \mathbbm{1}_{\mathcal{N}_k (x_n)} (X_n^j) 
	\end{align*}
	where $\mathcal{N}_k({x}_n)$ denotes the set of $k$ nearest neighbors on $x_n$, given by 
	\begin{align*}
		\mathcal{N}_k({x}_n)= \{ (X_n)_{i=1}^S : \sum_{j=1}^S \mathbbm{1}[\lVert {x}_n-X_n^i \rVert \geq \lVert x_n-X_n^j \rVert] \leq k \}
	\end{align*}	
	In conventional P2P energy trading scheme, prosumers negotiate their estimated energy deficits or surpluses with their neighbors to determine the optimal P2P energy exchange, regardless of the biased or asymmetrical estimates or contextual data used in their decision-making. To improve decision-making in this distributed framework, a distributed kNN algorithm that incorporates data from neighboring prosumers and adjusts weights accordingly are proposed as follows
	\begin{align*}
		w_n^i(x_n) = &\gamma_n \sum_{j=1}^S \frac{1}{k} \mathbbm{1}_{\mathcal{N}_k (x_n)} (X_n^j) \nonumber \\ 
		&+ (1-\gamma_n) \sum_{j=1}^S \frac{1}{\lvert \mathcal{N}_n \rvert} \mathbbm{1}_{\mathcal{N}_k (\hat{x}_{mn})} (X_n^j)		
	\end{align*}
	Here, $0 \leq \gamma_n \leq 1$ represents the degree of reliance on a prosumer's own covariate data relative to that of its neighbors.

	\subsection{Distributed Implementation}
	The optimization problem \eqref{CSOKNN} is inherently a bi-level problem, challenging to solve. By substituting the lower-level problem with its Karush-Kuhn-Tucker (KKT) conditions, it can be reformulated into a single-level optimization problem. The complementarity conditions are linearized using bigM method \cite{fortuny1981representation}. To maintain autonomy and privacy in decision-making, a distributed consensus ADMM approach are implemented with an acceleration mechanism as suggested by \cite{goldstein2014fast}.
	
	Rewriting the decision variables as $z_n=\{p_{n}^b, p_{n}^s, p_n^{mt}, q_n^{mt}\}$ and $z_{nm}=p_{nm}$, we introduce auxiliary variables $\hat{z}_{nm} = p_{nm}$ to express the reciprocity constraints \eqref{recipe1} as
	\begin{subequations}\label{re}
	\begin{align}
		\hat{z}_{nm} + \hat{z}_{mn} &= 0, \quad\quad \forall m \in \mathcal{N}_n \label{re1} \\
		z_{nm} &= \hat{z}_{nm}, \quad \forall m \in \mathcal{N}_n \label{re2}
	\end{align}
	\end{subequations}
	The augmented Lagrangian for this problem, with $\lambda_{nm}$	as the dual variables associated with the reciprocity constraints \eqref{re2}, can be expressed as
	\begin{align} \label{CKNN}
		\mathcal{L}_n(z_n(X_n^i&,k), z_{nm}(X_n^i,k), \lambda_{nm}) = \nonumber \\
		& \sum_{i=1}^S \{C(z_n(X_n^i, k), z_{nm}(X_n^i, k); Y_n^i) \nonumber \\ 
		&+ (\lambda_{nm})^{\top} (z_{nm}-\hat{z}_{nm}) + \frac{{\rho}}{2}\|z_{nm}-\hat{z}_{nm}\|^2 \} 
	\end{align}
	
	Due to the decomposable nature of the augmented Lagrangian and constraint (15b), each prosumer $n$ can independently update its local variable $z_n$ in a fully decentralized manner by solving problem \eqref{zupdate}, while keeping the auxiliary variable $z_{nm}^\nu$ and Lagrange multiplier $\lambda_{nm}^\nu$ fixed.
	\begin{align} \label{zupdate}
		z_n = \arg\min_{z_n\in\mathcal{Z}_n} \mathcal{L}_n(z_n(X_n^i,k), z_{nm}^\nu (X_n^i,k), \lambda_{nm}^\nu) 
	\end{align}
	After solving \eqref{zupdate}, each prosumer $n$ broadcasts $z_{nm}^{\nu+1}$ to prosumer $m$.
	
	The auxiliary variables $\hat{z}_{nm}$ are updated by solving the problem \eqref{auxUpdate}, which is also decomposable, while keeping other variables fixed.
	\begin{align}\label{auxUpdate}
		\min_{\hat{z}_{nm}} -(\lambda_{nm}^\nu)^{\top}& \hat{z}_{nm} -(\lambda_{mn}^\nu)^{\top} \hat{z}_{mn} \nonumber \\ 
		+ \frac{\rho}{2} & \lVert \hat{z}_{nm} - z_{nm}^{\nu+1} \rVert_2^2 + \frac{\rho}{2} \lVert \hat{z}_{mn} - z_{mn}^{\nu+1} \rVert_2^2 \\
		\text{subject to}& \quad \eqref{re1} \nonumber
	\end{align}
	
	The update for the dual variables can be expressed as
	\begin{align}
		\lambda_{nm}^{i,\nu+1}=\lambda_{nm}^{i,\nu} + \rho(z_{nm}^{i,\nu+1}-\hat{z}_{nm}^{i,\nu+1}), \quad \forall m \in \mathcal{N}_n \label{DDA}
	\end{align}
	The detailed steps for each prosumer $n$ are presented in Algorithm \ref{alg1}.
	
	A detailed derivation of the KKT optimality conditions for the lower-level problem and the closed-form auxiliary variables update is provided in Appendix A.	
	By adopting this distributed strategy, prosumers only need to exchange P2P energy estimates and current observed covariate data with their neighbors. 
	
	\begin{algorithm}[t]
		\caption{Distributed Energy Trading Algorithm}\label{alg1}
		\begin{algorithmic}[1]
			\State \textbf{Initialize} \quad $\forall n\in\mathcal{N}$, $v=1$, $z_n^{i,1} \in \mathcal{Z}_n, \lambda_{n} \in \mathbb{R}^{N_n}$,
			\State Train a local KNN model using data points $(x_n^j,y_n^j)_{j=1}^S$,
			\Repeat \quad $\text{for each} \quad n\in\mathcal{N}$	
			\State \textbf{Primal Update}
			\State Send new data $\hat{x}_{nm} = x_{n}$ to all neighbors $m \in \mathcal{N}_n$,
			\State Update the decision variables $z_n^{\nu+1}$ by \eqref{zupdate},
			\State \textbf{Auxiliary Update}
			\State Update the auxiliary variables  $z_{nm}^{\nu+1}$ by \eqref{auxUpdate}
			\State Send $z_{nm}^{\nu+1}$ to all neighbors $m \in \mathcal{N}_n$,
			\State \textbf{Dual Update}
			\State Update the dual variables $\lambda_{nm}^{\nu+1}$ by \eqref{DDA} $\forall m \in \mathcal{N}_n$,
			\State \textbf{Acceleration Step}
			\State \hspace{\algorithmicindent} $a_{nm}^{\nu+1}=\frac{1+\sqrt{1+(4a_{nm}^{\nu})^2}}{2}$
			\State \hspace{\algorithmicindent} $\bar{\lambda}_{nm}^{\nu+1} = \lambda_{nm}^{\nu+1} + \frac{a_{nm}^{\nu}-1}{a_{nm}^{\nu+1}}(\lambda_{nm}^{\nu+1}-\lambda_{nm}^{\nu})$ 
			\State $\nu \gets \nu + 1$ 
			\Until{$\max(\lVert z_{nm}^{\nu+1}-\hat{z}_{nm}^{\nu} \rVert) \leq \epsilon_N$ }
		\end{algorithmic}
	\end{algorithm}

	\section{Numerical Results and Discussion}\label{sec4}
	This section presents the numerical results of our proposed distributed contextual stochastic optimization algorithm. Its performance is evaluated in optimizing energy consumption, mitigating peak loads, and enhancing grid stability.
	First, the simulation setup, including dataset, parameters, and implementation are described. Next, results of the case study are analyzed to demonstrate the algorithm's effectiveness. Finally, sensitivity analysis are conducted to compare our approach with existing methods.

	\subsection{Experimental Setup}
	The performance of our proposed method (Algorithm \ref{alg1}) was evaluated using a covariate dataset from Tehran, Iran, which comprises solar PV generation, residential/commercial consumption patterns collected from June 2020 to 2023, and weather related data including time-of-day, day-of-week, temperature, cloud cover, and solar irradiance were sourced from the Iranian national meteorological organization (INMO) \cite{irimo_dataset_2024}.
	
	Real-time and P2P electricity prices in Tehran, where direct historical price data is unavailable, a transfer learning approach is adapted using an auto-regressive integrated moving average with exogenous variables (ARIMAX). The ARIMAX model, was trained on normalized historical data from Los Angeles \cite{eia_electricity_market_data}\cite{visualcrossing_la_weather} using MATLAB’s arima function, incorporating relevant covariates. To account for the unique mid-day cooling demand in Tehran, an additional covariate was created to reflect this increased demand. This model was then adapted and used to forecast Tehran's real-time prices, inferring missing price data across regions. The predicted real-time price using this method is illustrated in Fig. \ref{fig6}.
	The dataset was divided into two parts for training and testing purposes. For further analysis, data from previous years for the same month are employed. Shiftable loads were conservatively estimated to constitute $10-50\%$ of total consumption. 

	The study focuses on a typical energy community consisting of ten prosumers, with a 24-hour time horizon and a time step of $\Delta t=1$. Table I outlines the typical parameter values used in the model. The ADMM parameters were set to $\rho=0.5$ and $\epsilon_N=10^{-4}$.
	Algorithms \ref{alg1} and the optimization problems \eqref{CKNN} were implemented in MATLAB and solved using Mosek. All simulations were conducted on a desktop PC equipped with an Intel\textsuperscript{\textregistered} Core i7-1165G7 four-core CPU processor running at 4.70 GHz with 16GB of RAM.

\begin{figure}[!t] 
	\centering
	\subfloat[\label{1c}]{%
		\includegraphics[width=0.5\linewidth]{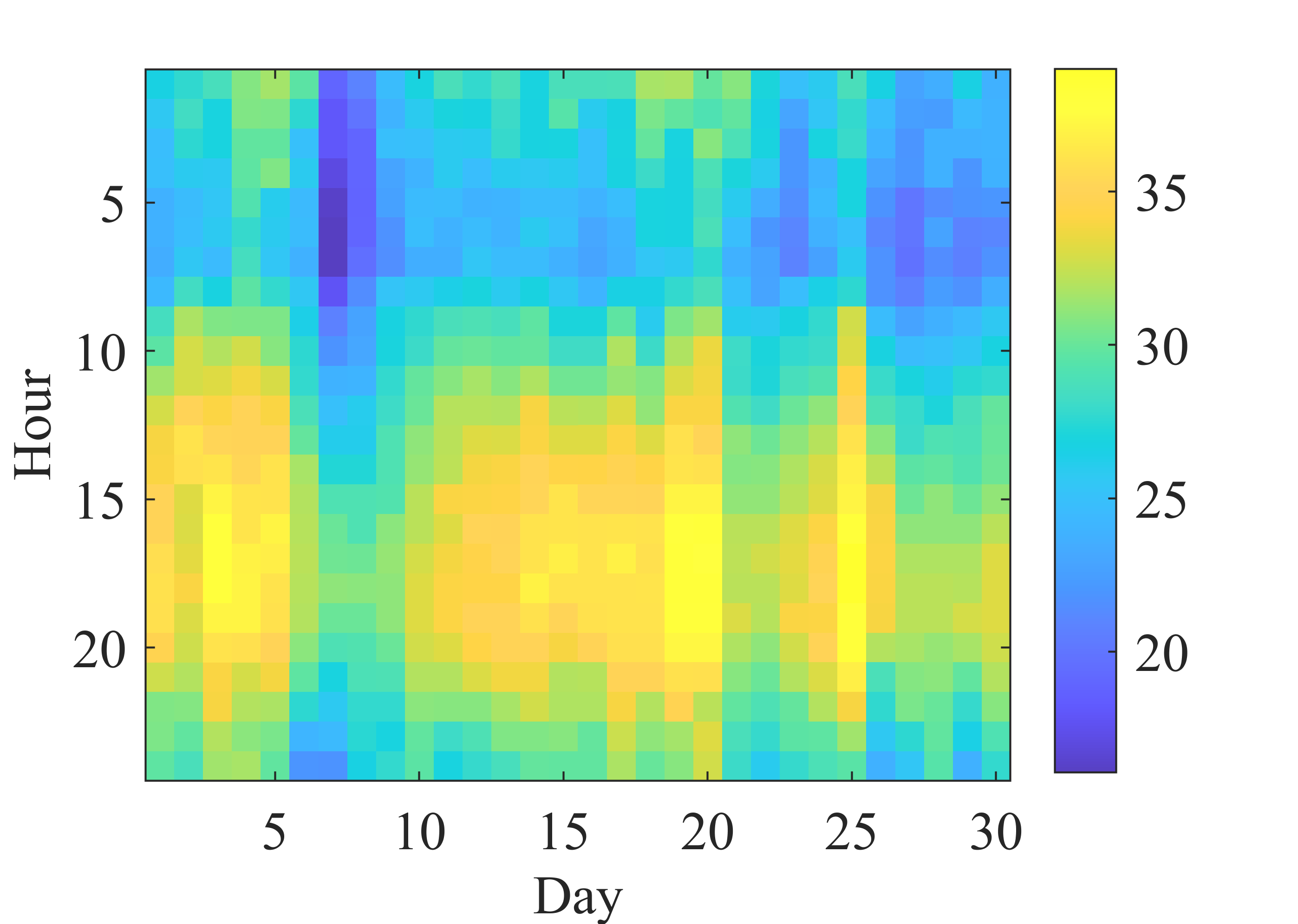}}
	\hfill
	\subfloat[\label{1d}]{%
		\includegraphics[width=0.5\linewidth]{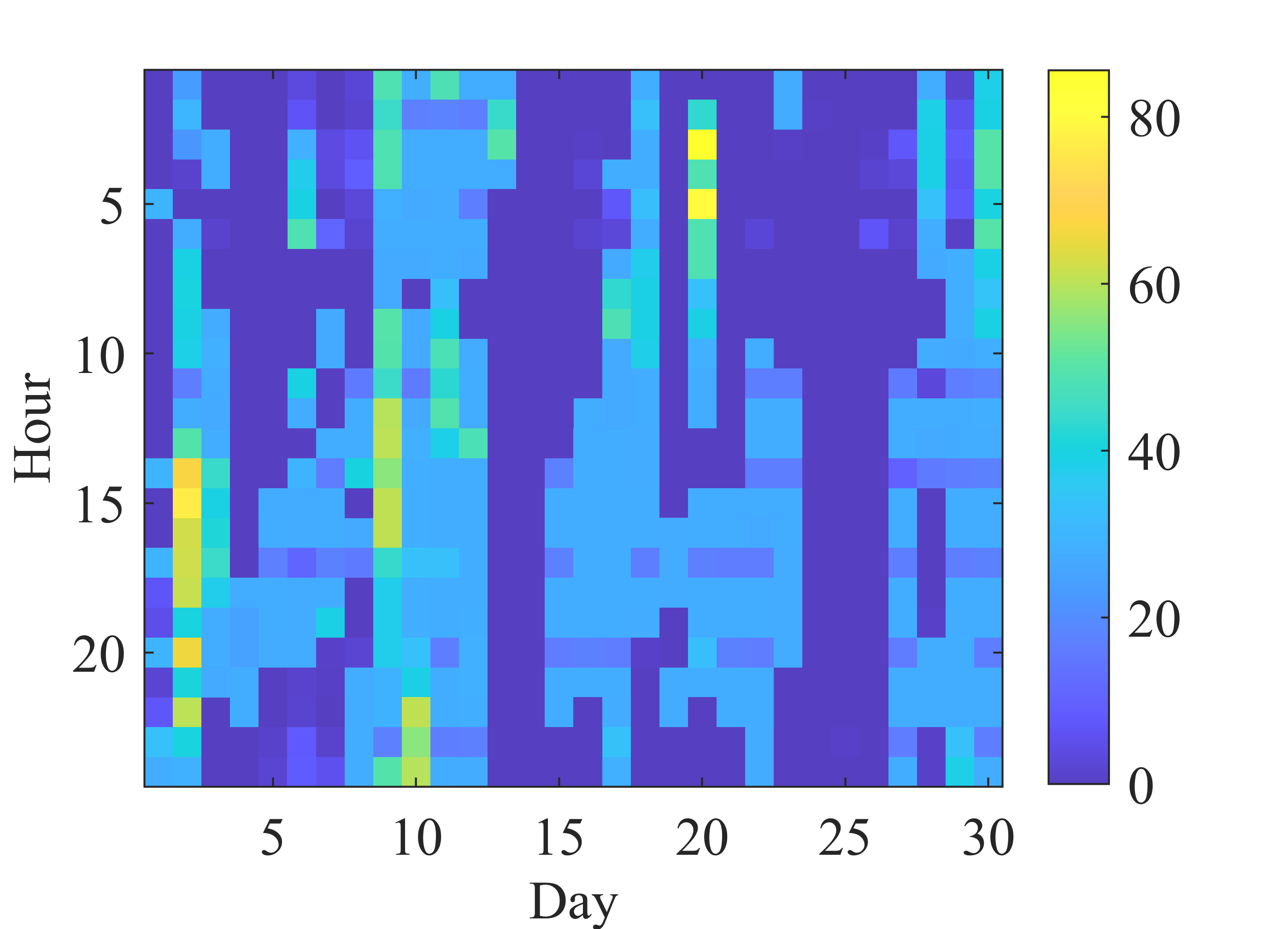}}
	\caption{(a) Temperature Data ($^\circ$C), (b) Cloud covering data ($\%$).}
	\label{fig3} 
\end{figure}

\begin{figure}[!t] 
	\centering
	\subfloat[\label{1a}]{%
		\includegraphics[width=0.5\linewidth]{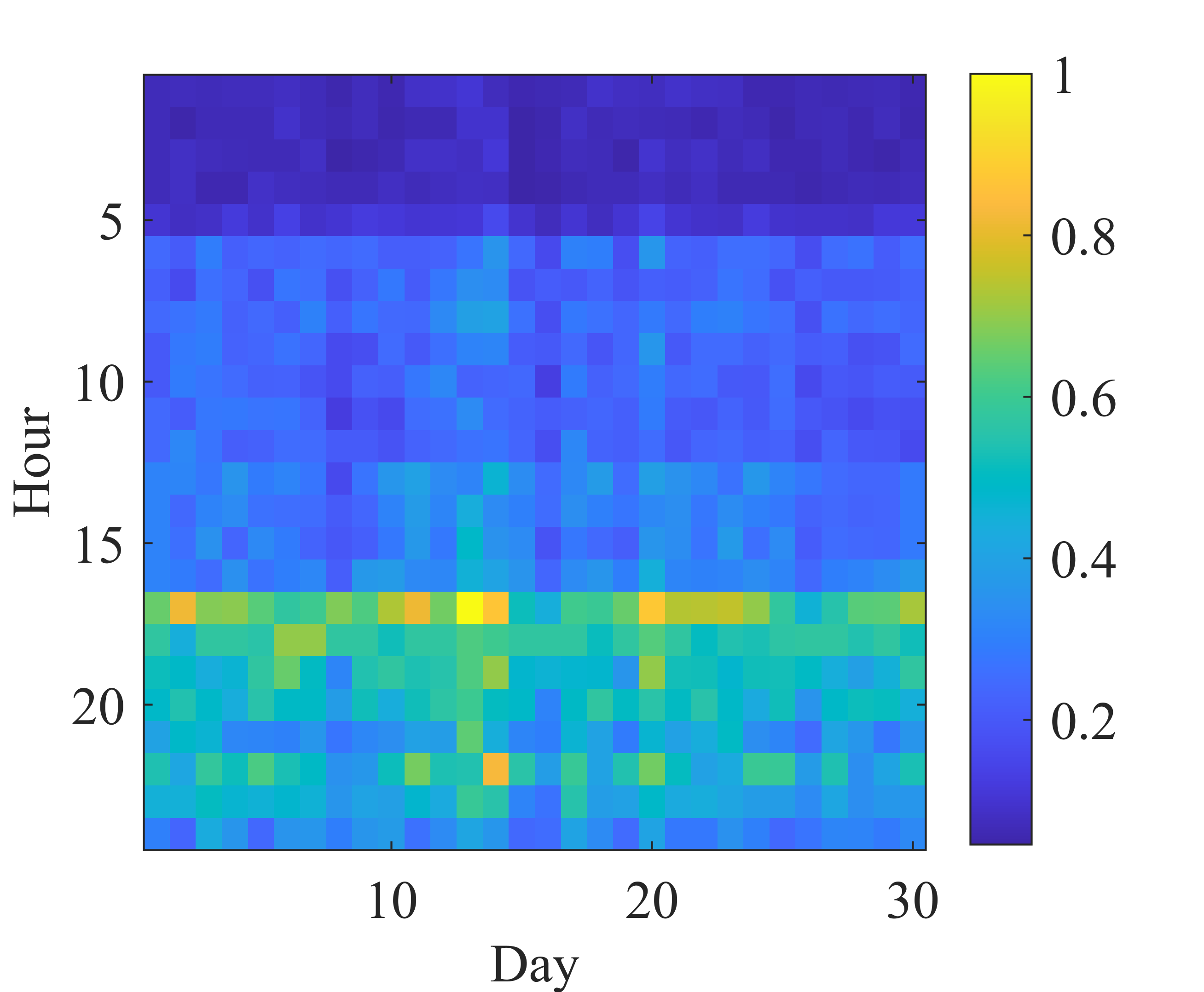}}
	\hfill
	\subfloat[\label{1b}]{%
		\includegraphics[width=0.5\linewidth]{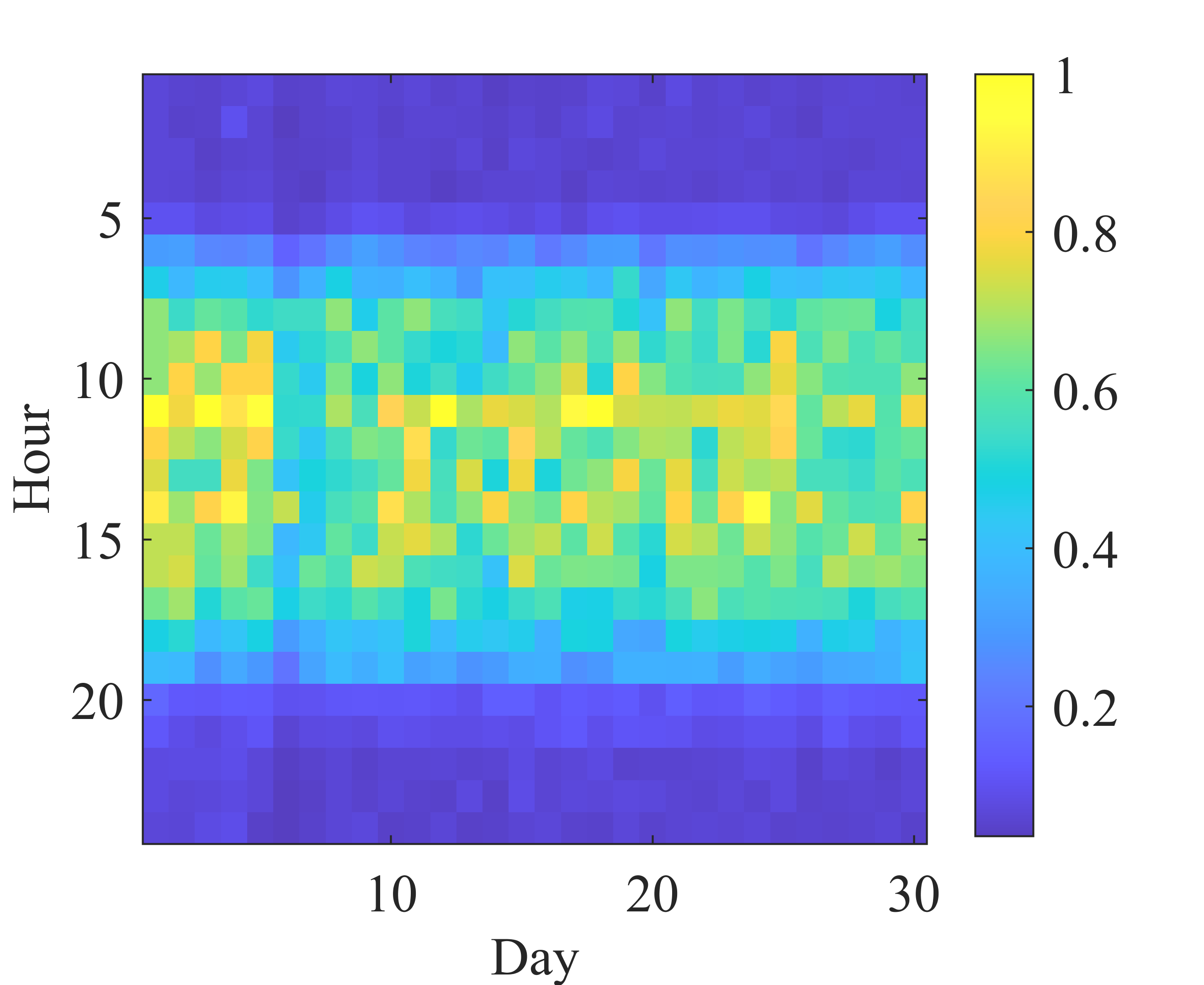}}
	\caption{(a) Normalized consumption (p.u.) for (a) residential Prosumer 1 and (b) commercial Prosumer 2.}
	\label{fig4} 
\end{figure}

\begin{figure}[!t]
	\centering
	\subfloat[\label{1a}]{%
		\includegraphics[width=0.5\linewidth]{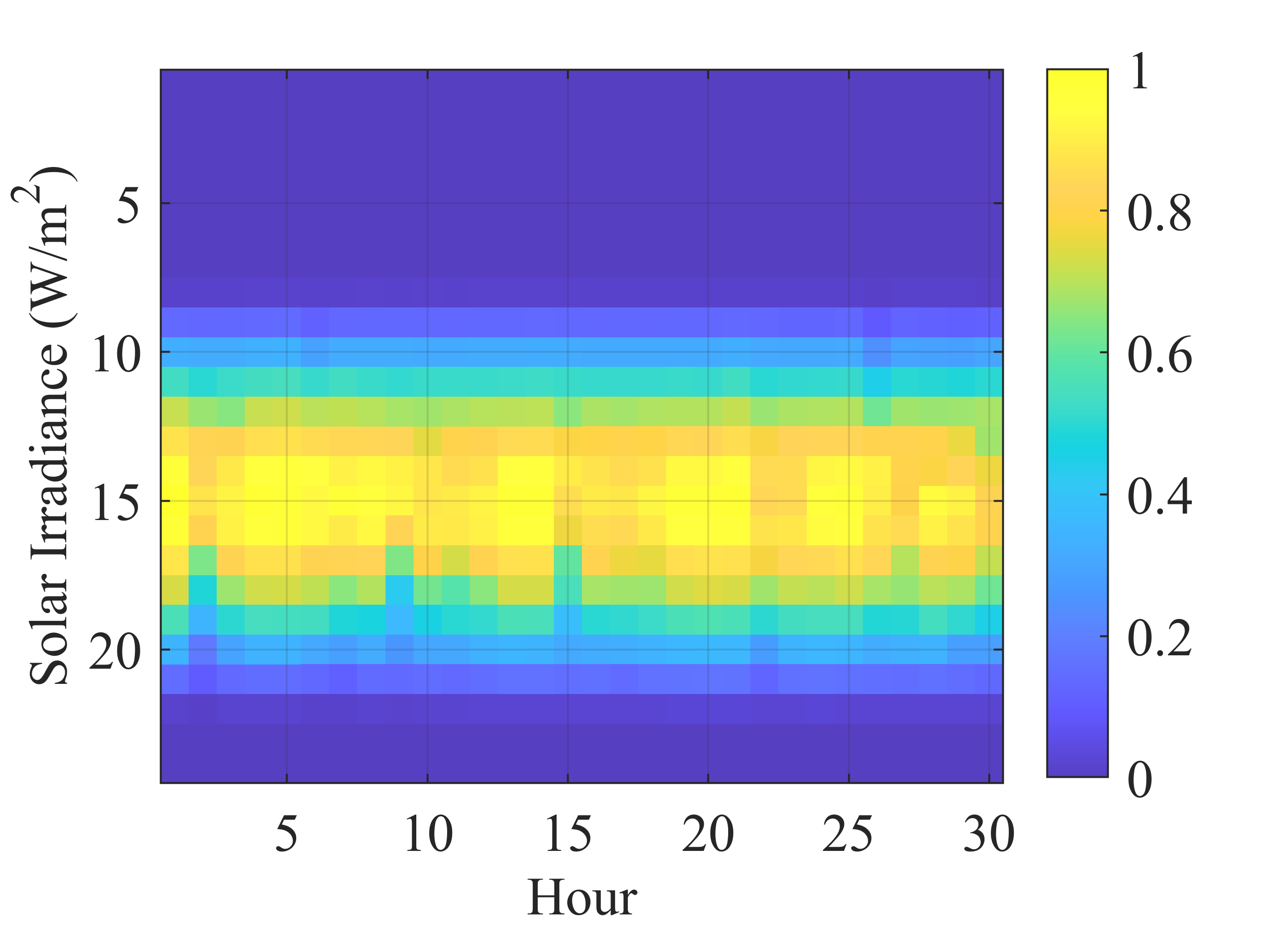}}
	\hfill
	\subfloat[\label{1b}]{%
		\includegraphics[width=0.5\linewidth]{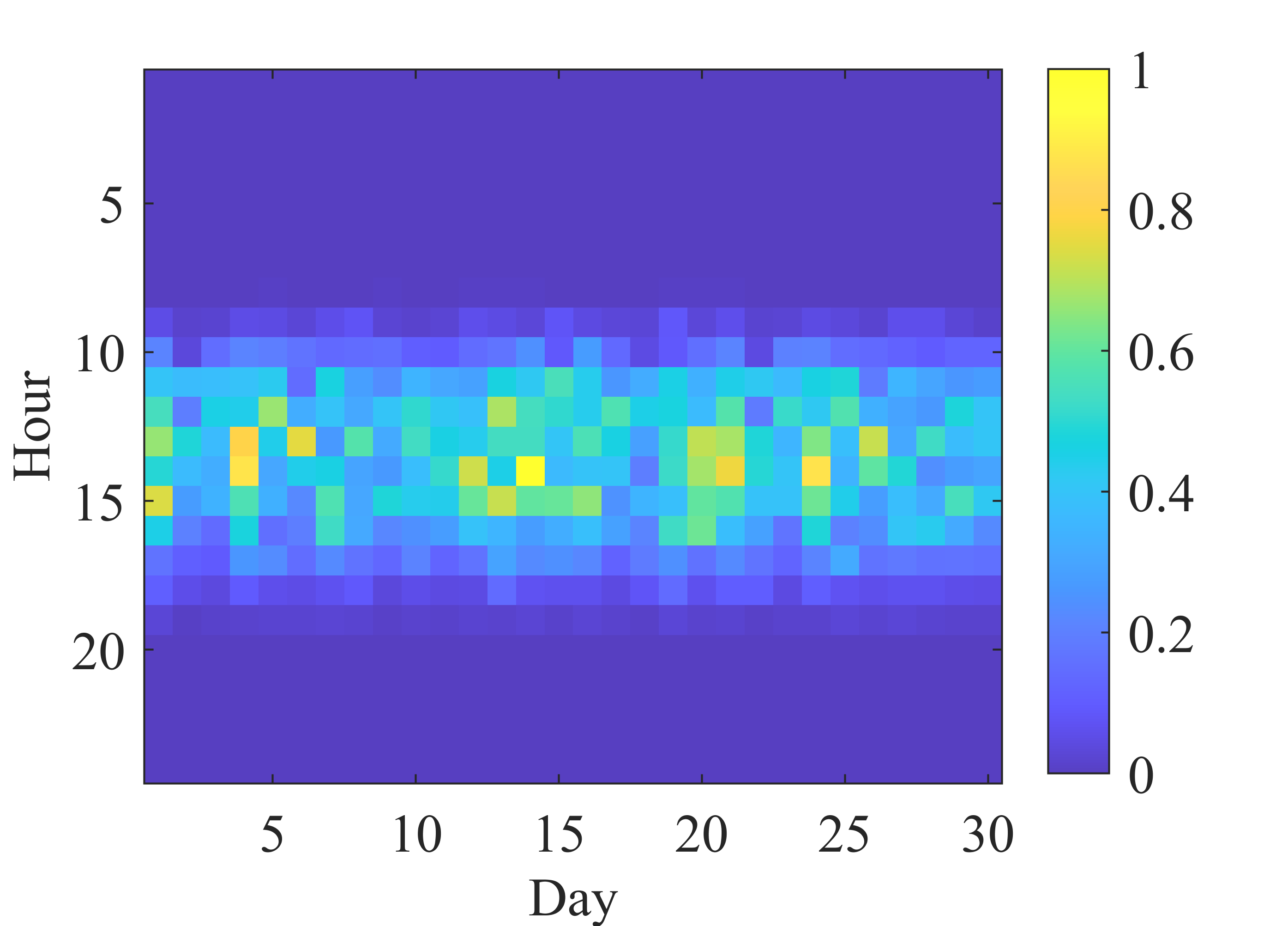}}
	\caption{(a) Normalized PV generation (kWh), and (b) solar irradiance (W/m$^2$).}
	\label{fig5} 
\end{figure}

\begin{figure}[!t]
	\centering
	\subfloat[\label{1a}]{%
		\includegraphics[width=0.5\linewidth]{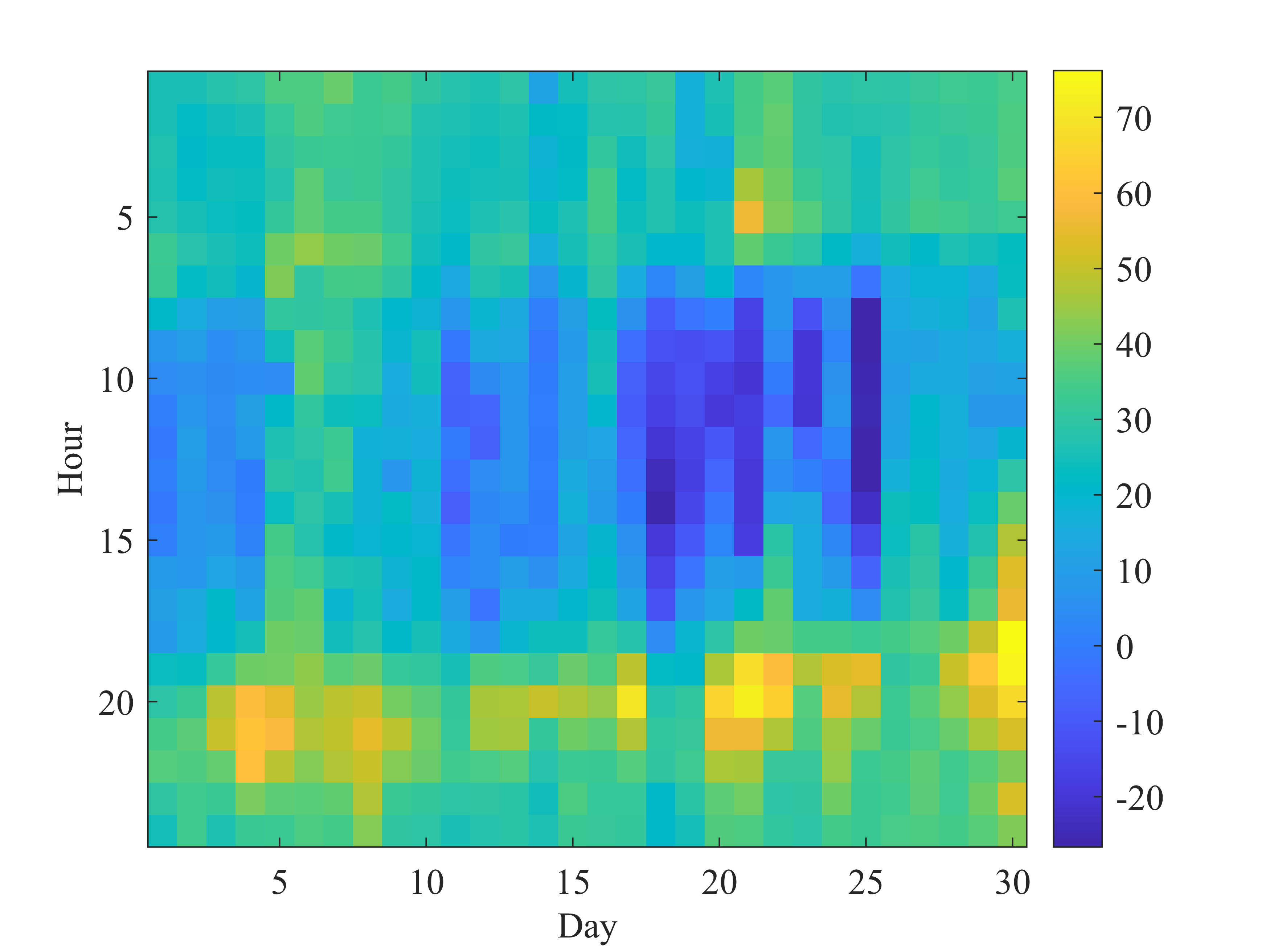}}
	\hfill
	\subfloat[\label{1b}]{%
		\includegraphics[width=0.5\linewidth]{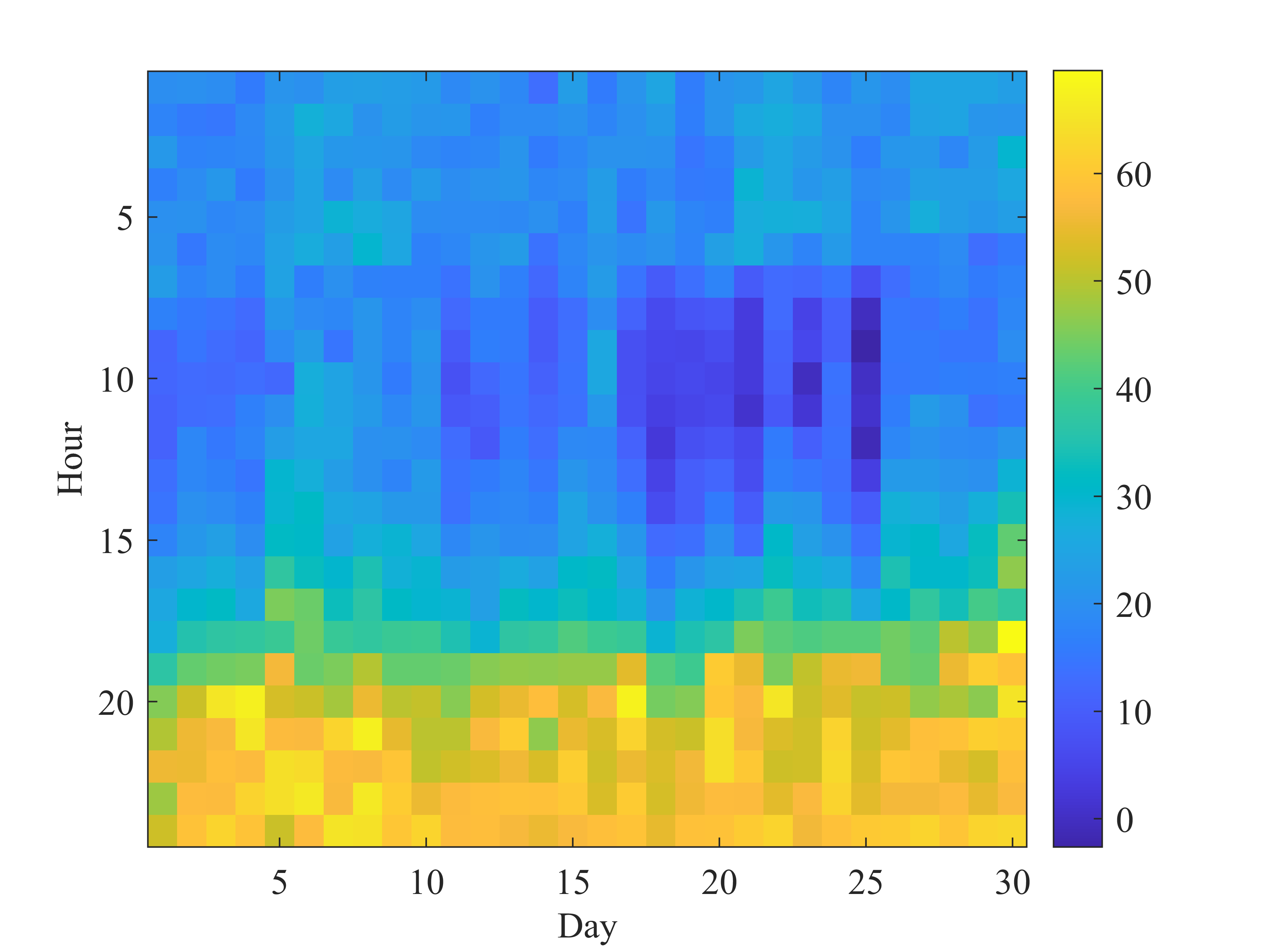}}
	\caption{(a) Actual price (LA), and (b) Forecast price (TEH).}
	\label{fig6} 
\end{figure}

	\begin{table}[!htb]
		\centering
		\caption{Typical Parameter Values}\label{tab1}%
		\begin{tabular}{@{}llll@{}}
			\hline\hline 
			Parameter & Value  & Parameter & Value\\
			\hline 
			
			$\eta_n^c$/$\eta_n^d$ & $0.9/0.9$ & $\underline{e}/\overline{e}$ & 10/200 $kWh$  \\
			$\alpha_1/\alpha_2$ & $0.015/0.05 \$/kWh$ & $C^s$ & 150 $kWh$  \\
			$\alpha^b/\alpha^s$ & 0.1 \textcent$/kW^2$ & $\overline{p}^{e}$ & 5 $kW$  \\
			$\beta^b/\beta^s$ & 1 \textcent$/kW$ & $\overline{p}^s$ & 20 $kW$ \\
			$c^{e}$ & $[3-9]$ \textcent$/kW$ & $\overline{p}^{b}$ & 20 $kW$ \\		
			\hline
		\end{tabular}
	\end{table}

\begin{figure}[!t]
	\centering
	\subfloat[\label{1a}]{%
		\includegraphics[width=0.5\linewidth]{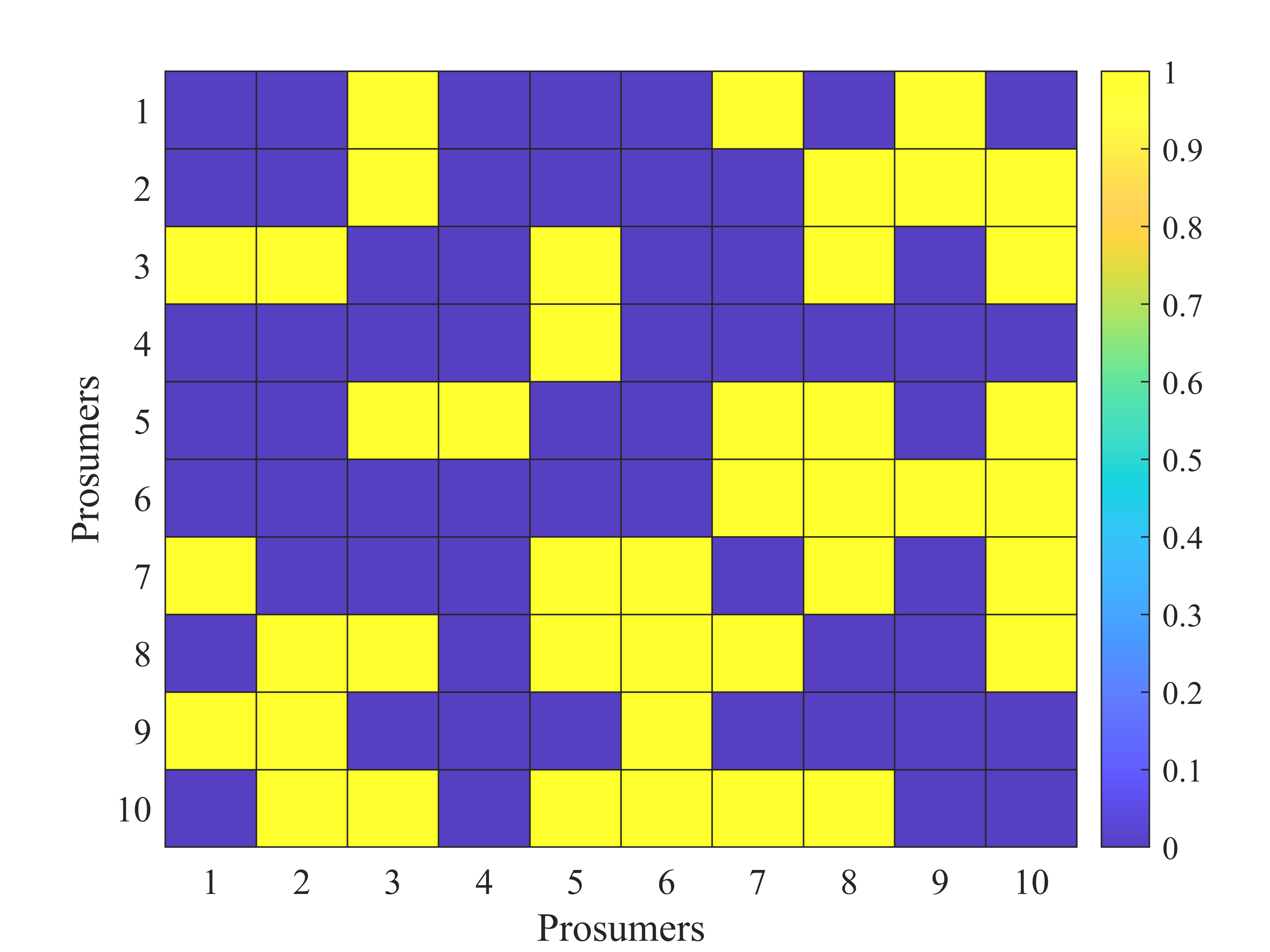}}
	\hfill
	\subfloat[\label{1b}]{%
		\includegraphics[width=0.5\linewidth]{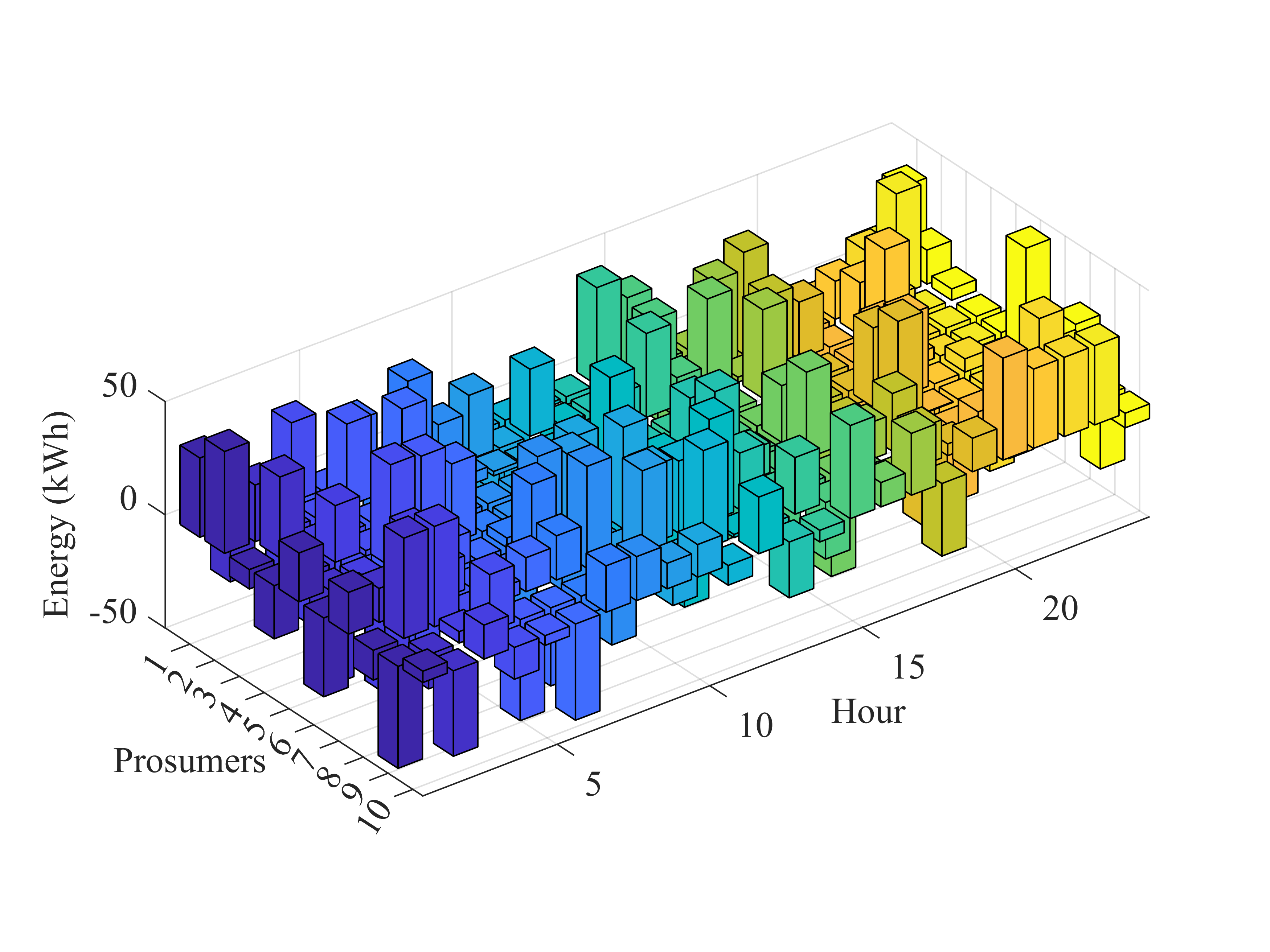}}
	\caption{EnergyExchange (a) Adjacency matrix of prosumers, (b) P2P energy transactions between prosumers.}
	\label{fig7} 
\end{figure}

\begin{figure}[!t]
	\centering
	\subfloat[\label{1a}]{%
		\includegraphics[width=0.5\linewidth]{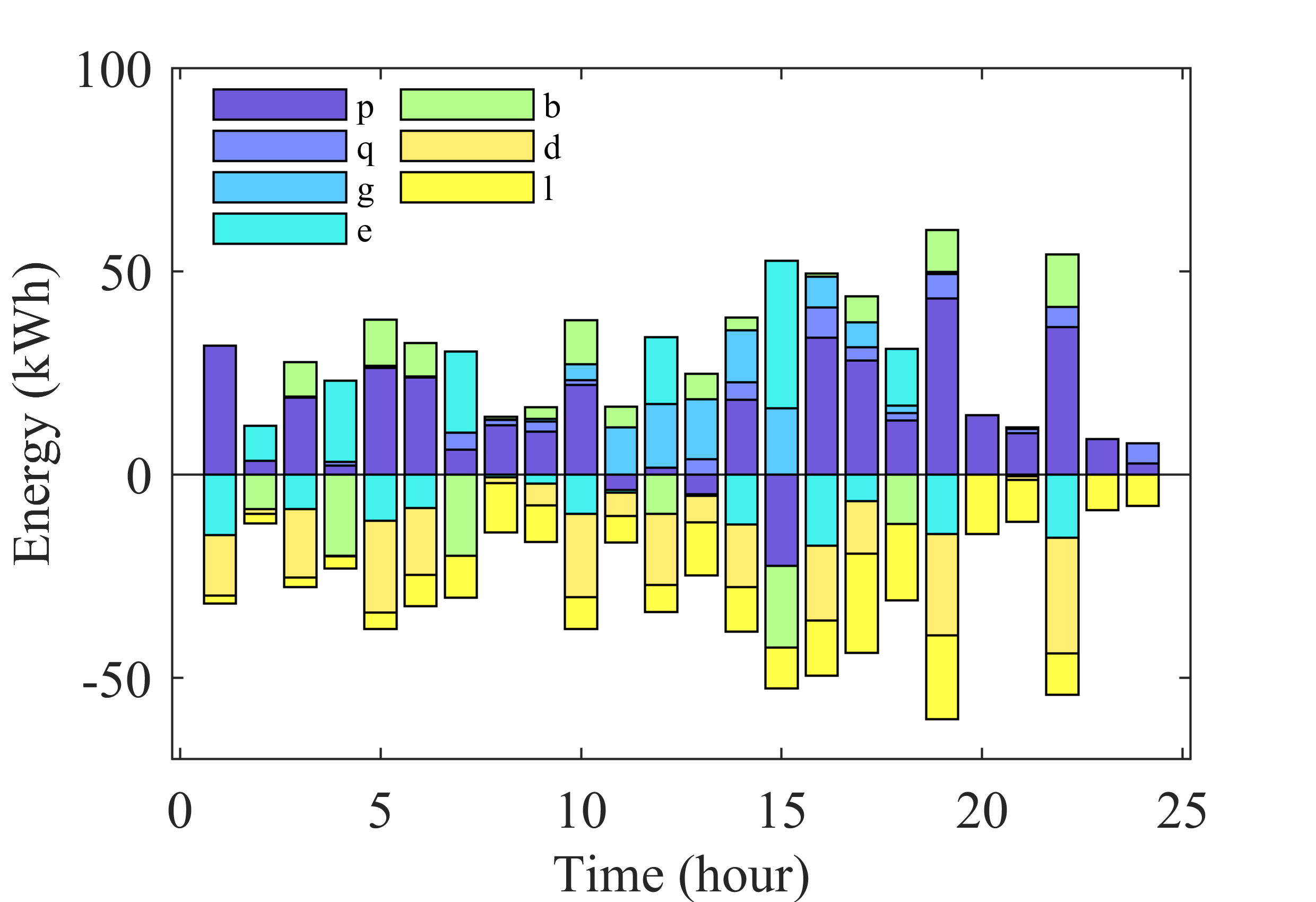}}
	\hfill
	\subfloat[\label{1b}]{%
		\includegraphics[width=0.5\linewidth]{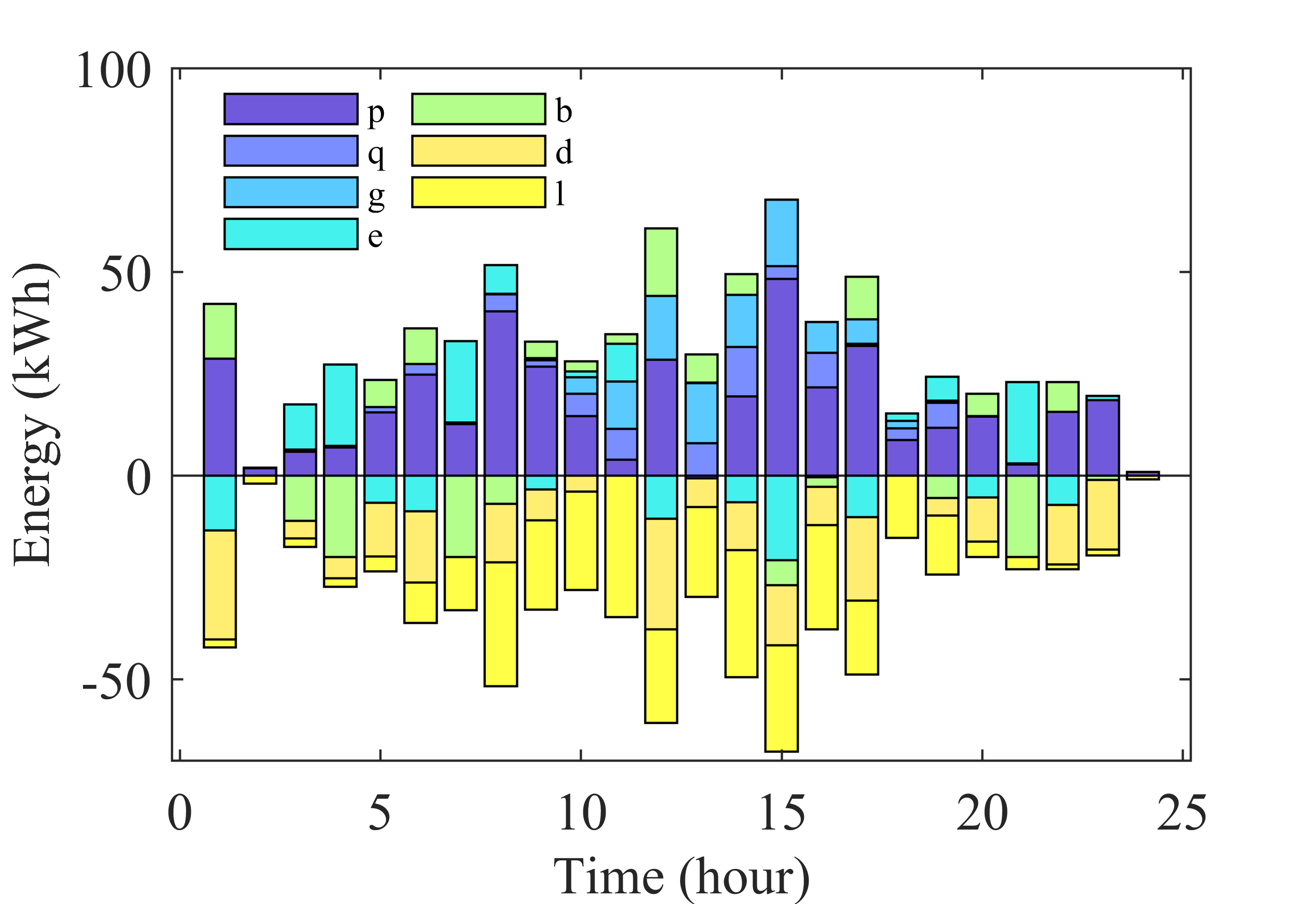}}
	\caption{Energy Management of Representative Prosumers (a) Prosumer 1 (Residential), and (b) Prosumer 2 (Commercial).}
	\label{fig8} 
\end{figure}

\begin{figure}[!t]
	\centering
	\subfloat[\label{1a}]{%
		\includegraphics[width=0.5\linewidth]{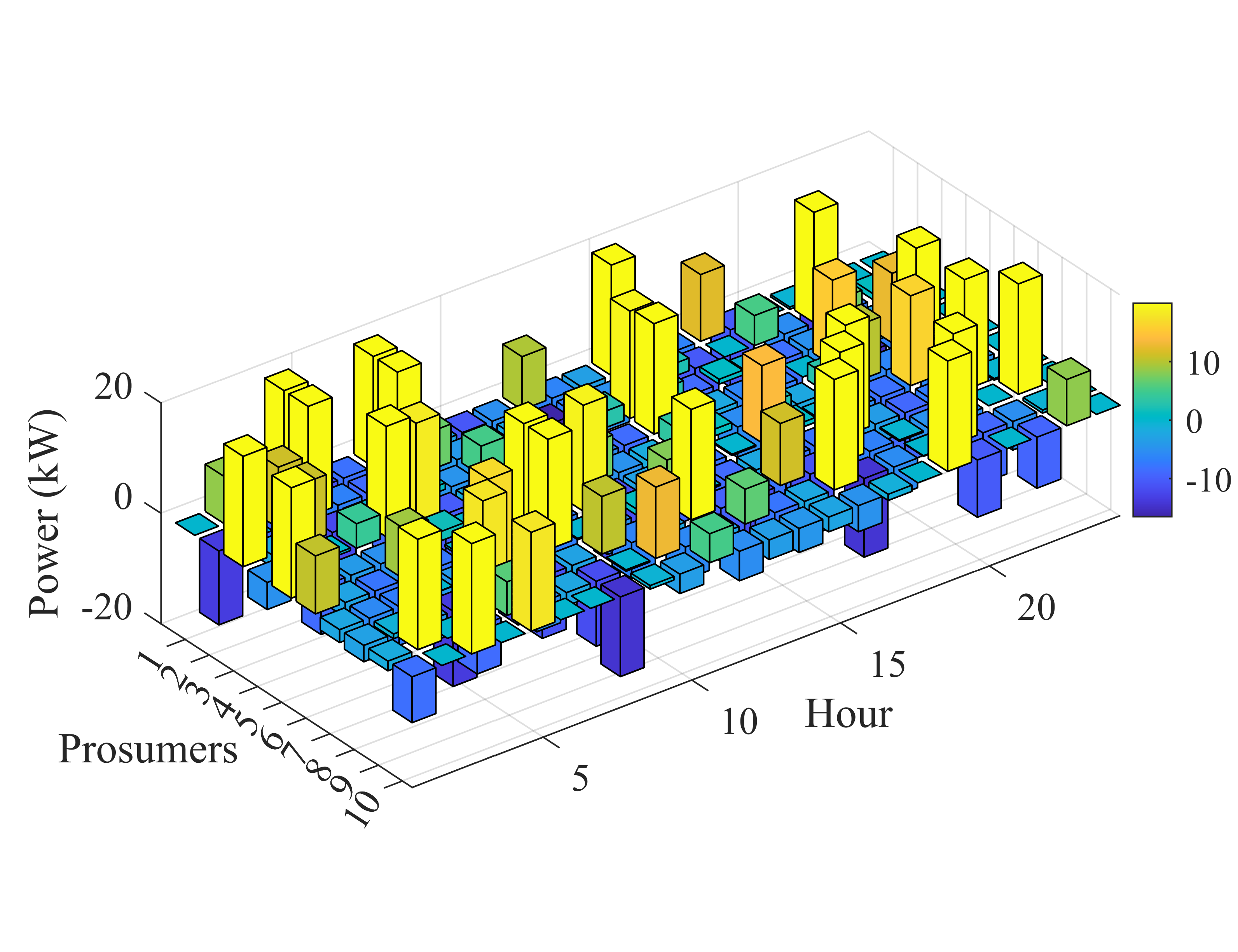}}
	\hfill
	\subfloat[\label{1b}]{%
		\includegraphics[width=0.5\linewidth]{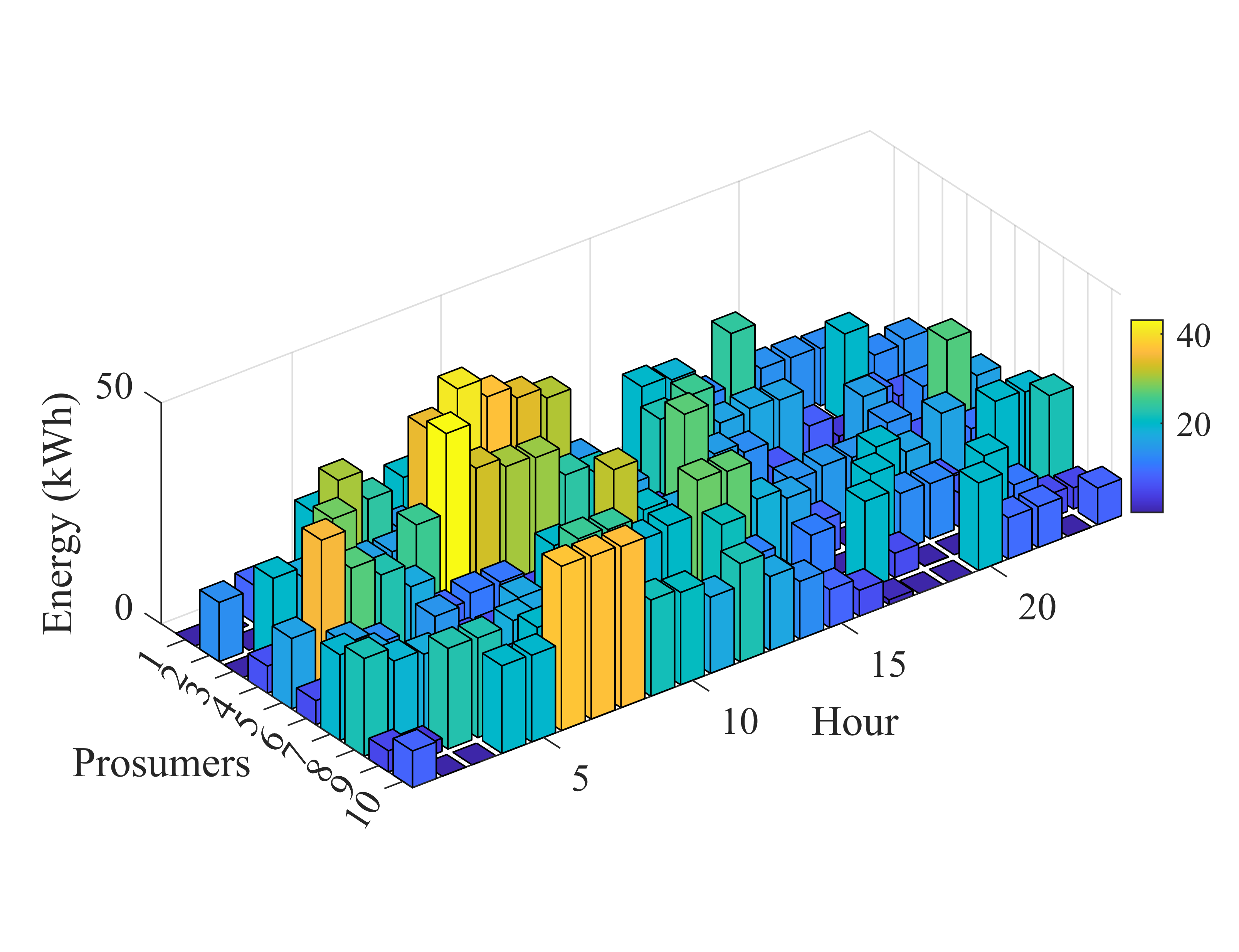}}
	\caption{Operation of BESSs (a) Power (kW), (b) State-of-charge (kWh).}
	\label{fig9} 
\end{figure}
\begin{figure}[!t]
	\centering
	\subfloat[\label{1a}]{%
		\includegraphics[width=0.5\linewidth]{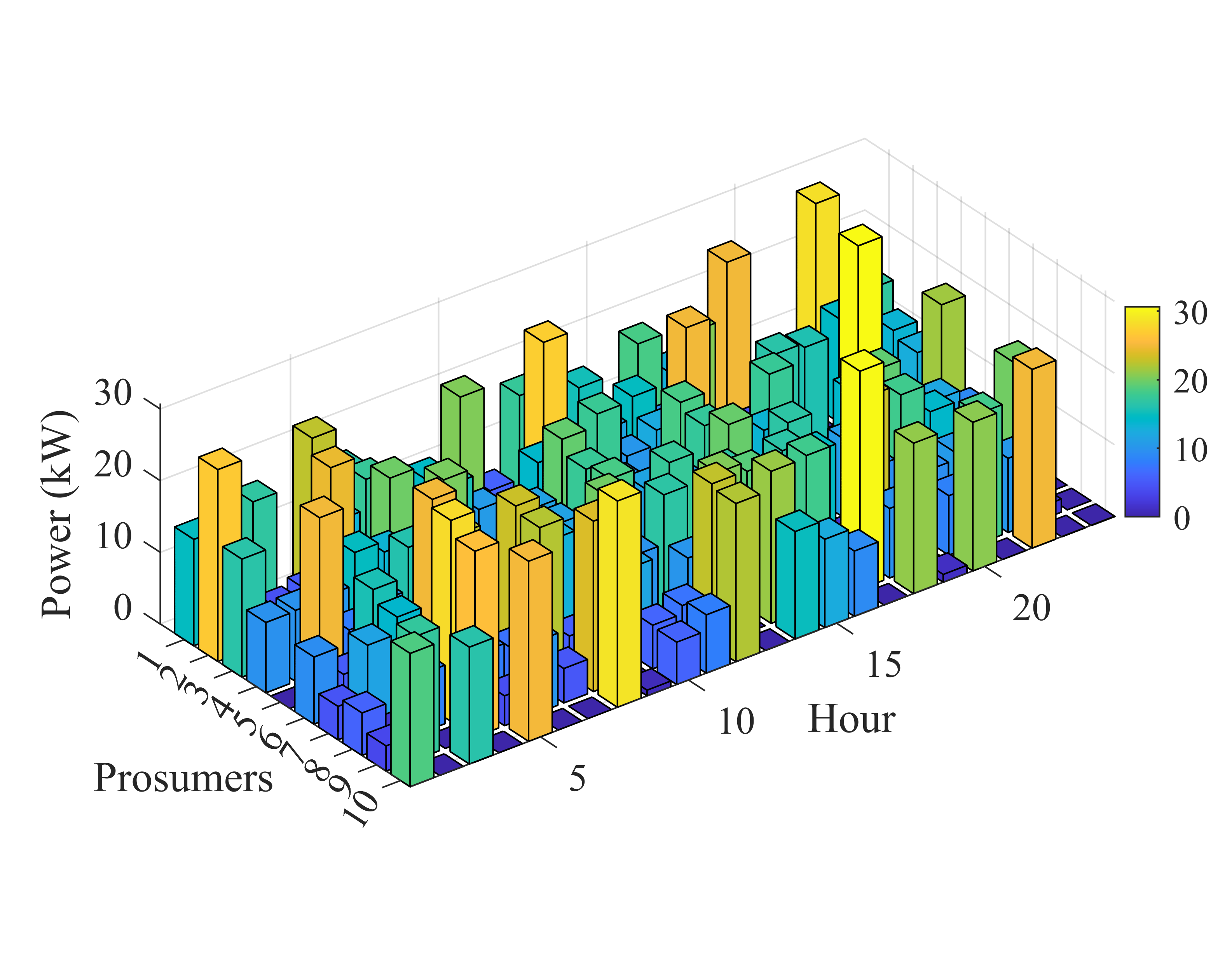}}
	\hfill
	\subfloat[\label{1b}]{%
		\includegraphics[width=0.5\linewidth]{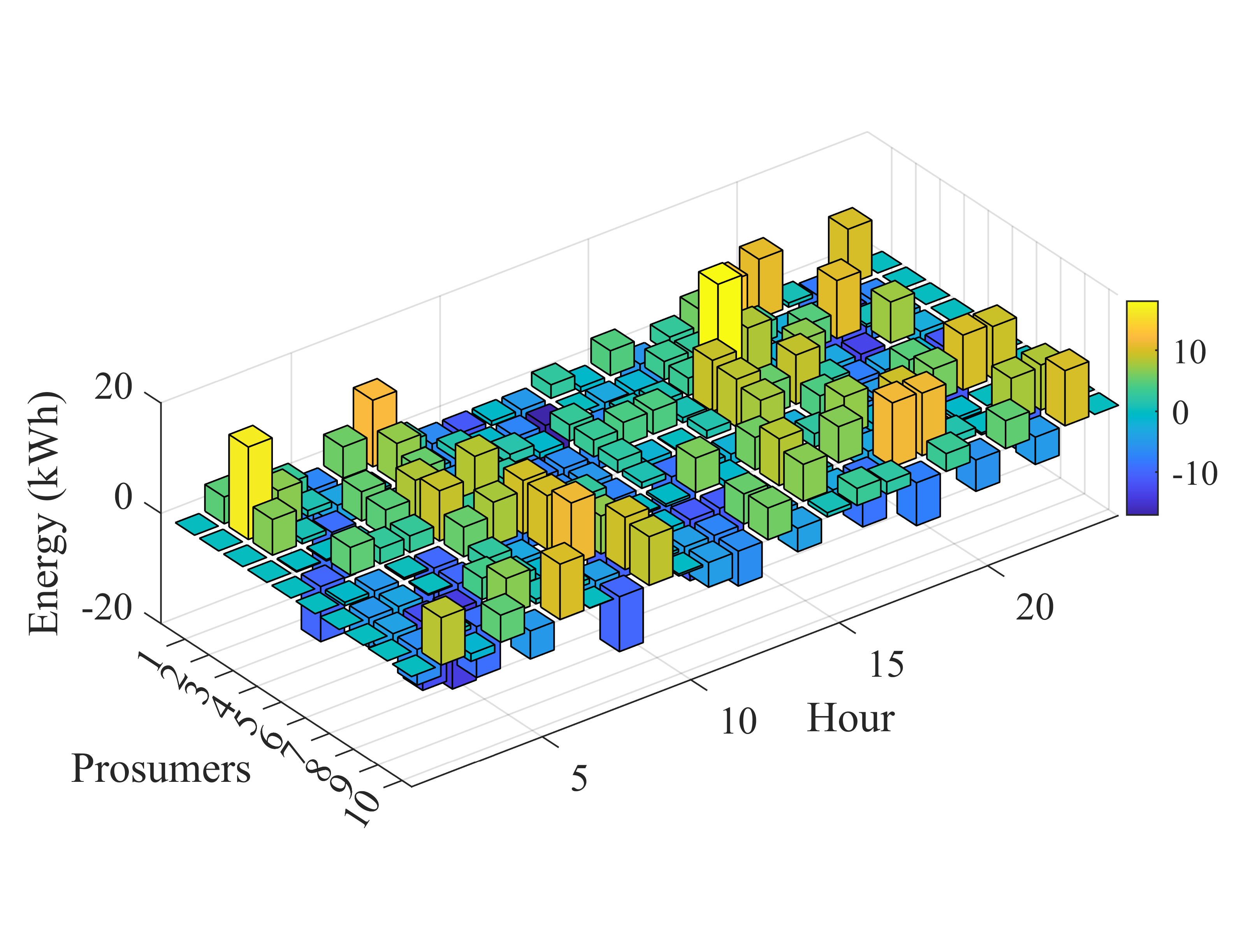}}
	\caption{Operation of SLs (a) Power (kW), (b) Energy shift state (kWh).}
	\label{fig10} 
\end{figure}
\begin{figure}[!t]
	\centering
	\subfloat[\label{1a}]{%
		\includegraphics[width=0.5\linewidth]{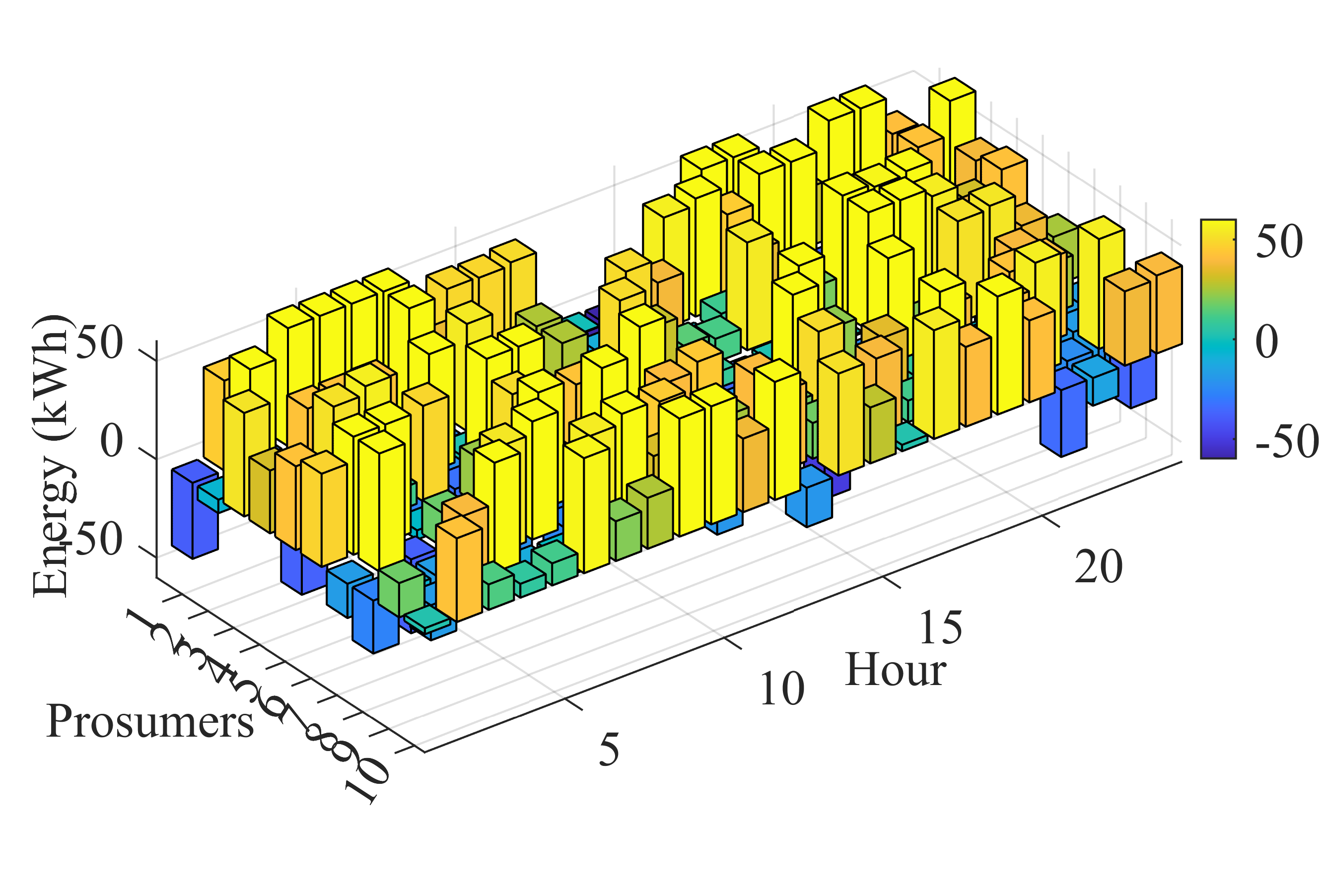}}
	\hfill
	\subfloat[\label{1b}]{%
		\includegraphics[width=0.5\linewidth]{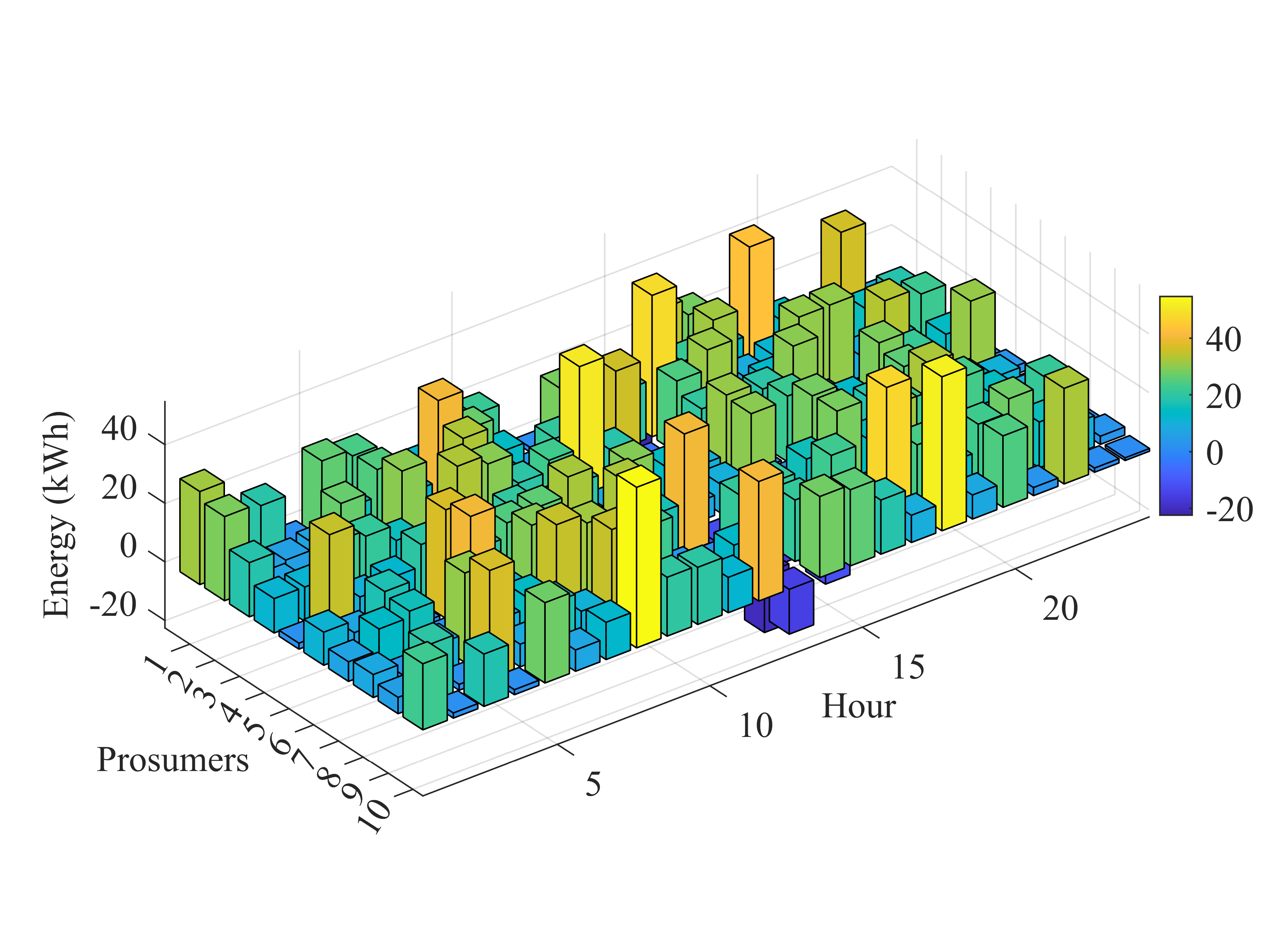}}
	\caption{Power exchange of prosumers with the main grid (a) Baseline scenario, (b) DCSO strategy}
	\label{fig101} 
\end{figure}

	\subsection{Results}
	
	The proposed algorithm was assessed through a case study involving data sharing and P2P energy trading within a defined network topology (Figure \ref{fig7}a). The resulting P2P energy transactions are depicted in Figure \ref{fig7}b. Prosumers optimized their flexible resources and trading strategies based on their own covariate data, as well as information exchanged with neighboring prosumers. The power profiles for two representative prosumers are presented in Figure \ref{fig8}.
	
	The effective utilization of flexible resources is essential for consistent energy trading in P2P as well as with the main grid scheme. Figures \ref{fig9}a and \ref{fig9}b present the power and state of charge of battery energy storage systems, while Figures \ref{fig10}a and \ref{fig10}b depict the power and energy shifting states of smart loads. The coordinated operation of these flexible resources significantly impacts power exchanges with the main grid. Figures \ref{fig101}a and \ref{fig101}b compare the power traded with the main grid in the baseline scenario and under the proposed framework.

	\subsection{Discussion}

	\subsubsection{Peak Load Mitigation}
\begin{figure}[!t]
	\centering
	\subfloat[\label{1a}]{%
		\includegraphics[width=0.5\linewidth]{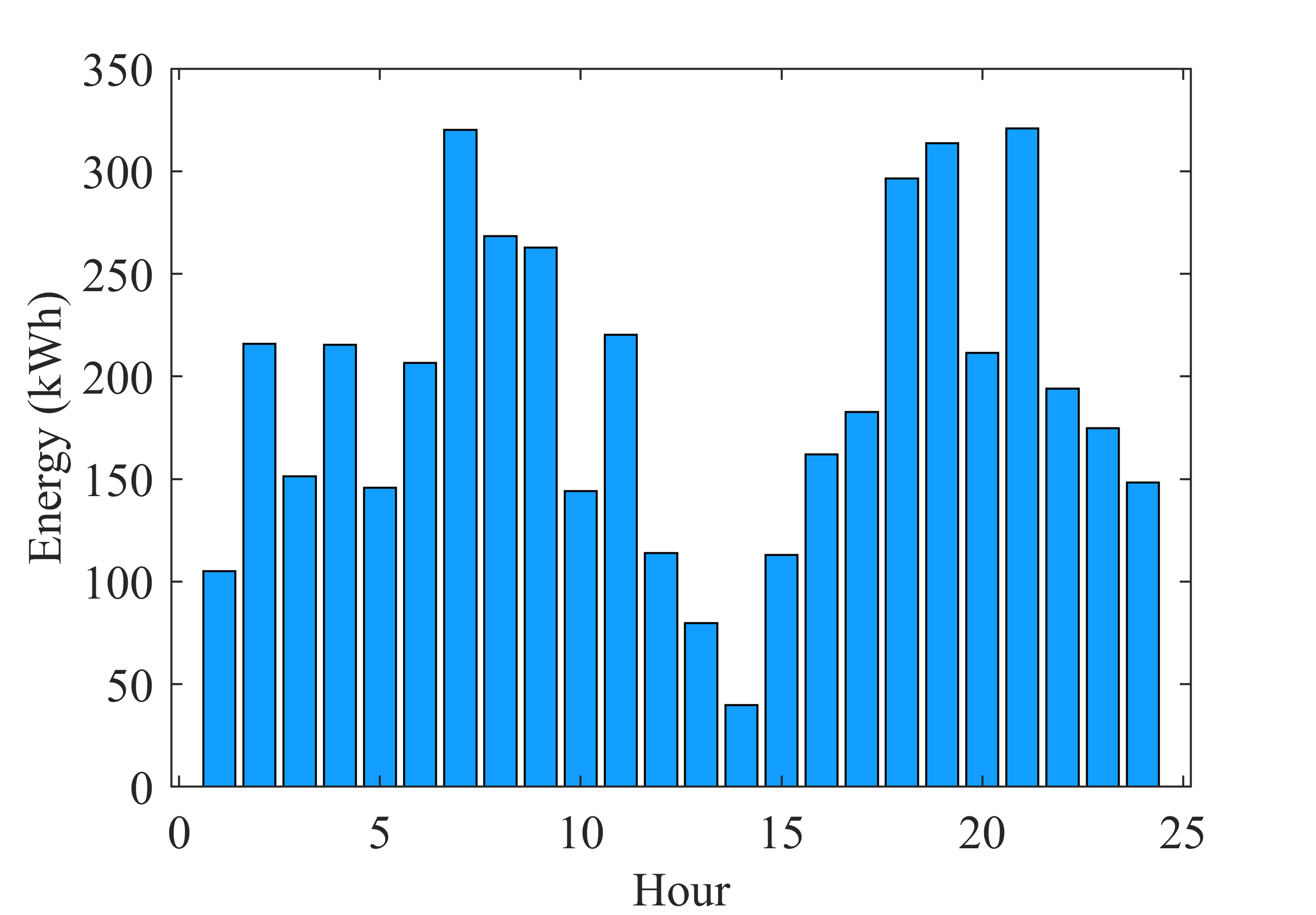}}
	\hfill
	\subfloat[\label{1b}]{%
		\includegraphics[width=0.5\linewidth]{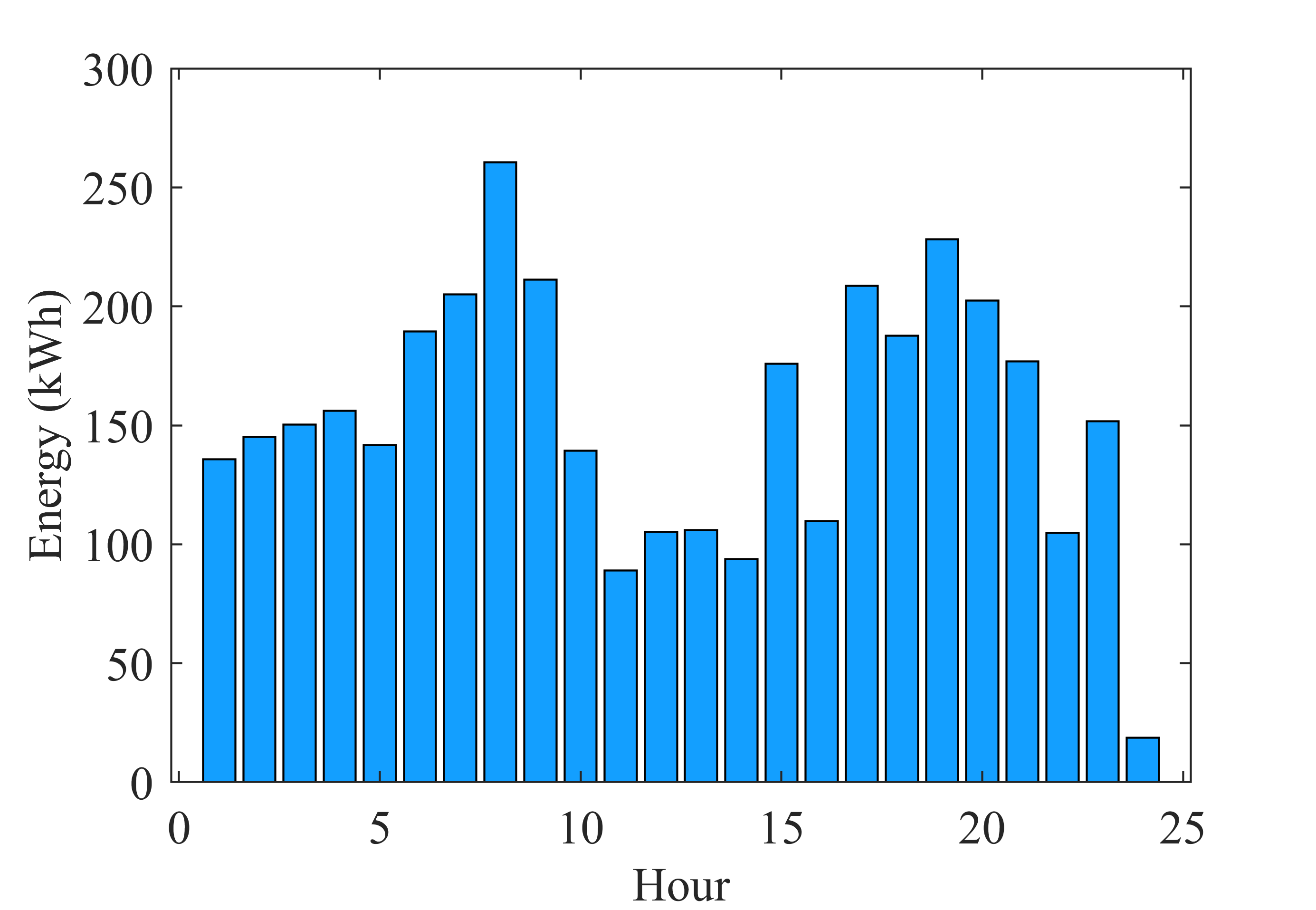}}
	\caption{Total power exchange with the main grid (a) Baseline scenario, (b) DCSO strategy}
	\label{fig102} 
\end{figure}

	Figure \ref{fig101}a illustrates that, without coordination, prosumers exhibit peak allowable load usage consistently throughout the day. This simultaneous demand from the main grid results in a significantly high total peak load, as shown in Figure \ref{fig102}a. However, the integration of P2P energy exchange and the strategic utilization of flexible resources leads to a significant reduction in grid load, with some scenarios achieving an average decrease of approximately 18\%. 

	\subsubsection{Economic Benefits}
	
	Table \ref{tab2} compares the proposed method (DCSO) with predict-then-optimize (PO), SAA, and SAA-KNN, highlighting its superior performance in handling uncertainties within the optimization process. DCSO consistently outperforms these baseline methods by up to 19\% and 7\%, respectively.
	Traditional accuracy-based methods often rely solely on predictions to determine decision variables, neglecting contextual information, which can lead to suboptimal outcomes. Our multi-time scale strategy effectively prescribes flexible resources and ensures consistent energy trading, contributing to power balance and smooth network operation.

	
	\begin{table}[!htb]
		\centering
		\caption{Average performance (N=50, I=50)}
		\label{tab2}
		\begin{tabular}{c c c c c}
			\hline\hline
			
			Method  & PO & SAA  & WSAA (KNN) & WSAA (CKNN) \\
			\hline
			
			Total cost (\textcent) & $11359$ & $10469$ & $9831$ & $9187$ \\			
			\hline
		\end{tabular}
	\end{table}

	\subsubsection{Sensitivity Analysis and Scalability}
		
\begin{figure}[!t]
		\centering
		\subfloat[\label{1a}]{%
			\includegraphics[width=0.5\linewidth]{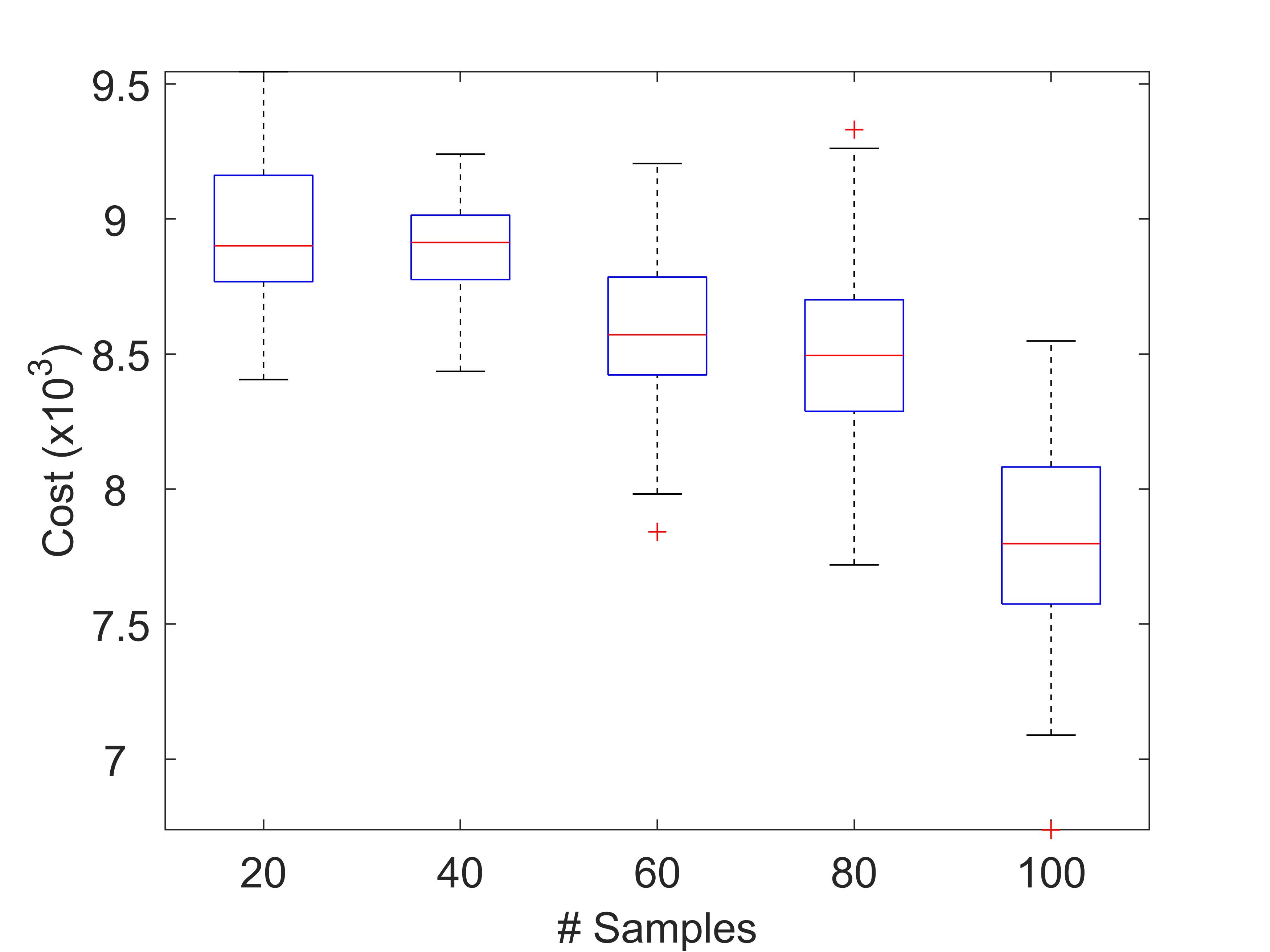}}
		\hfill
		\subfloat[\label{1b}]{%
			\includegraphics[width=0.5\linewidth]{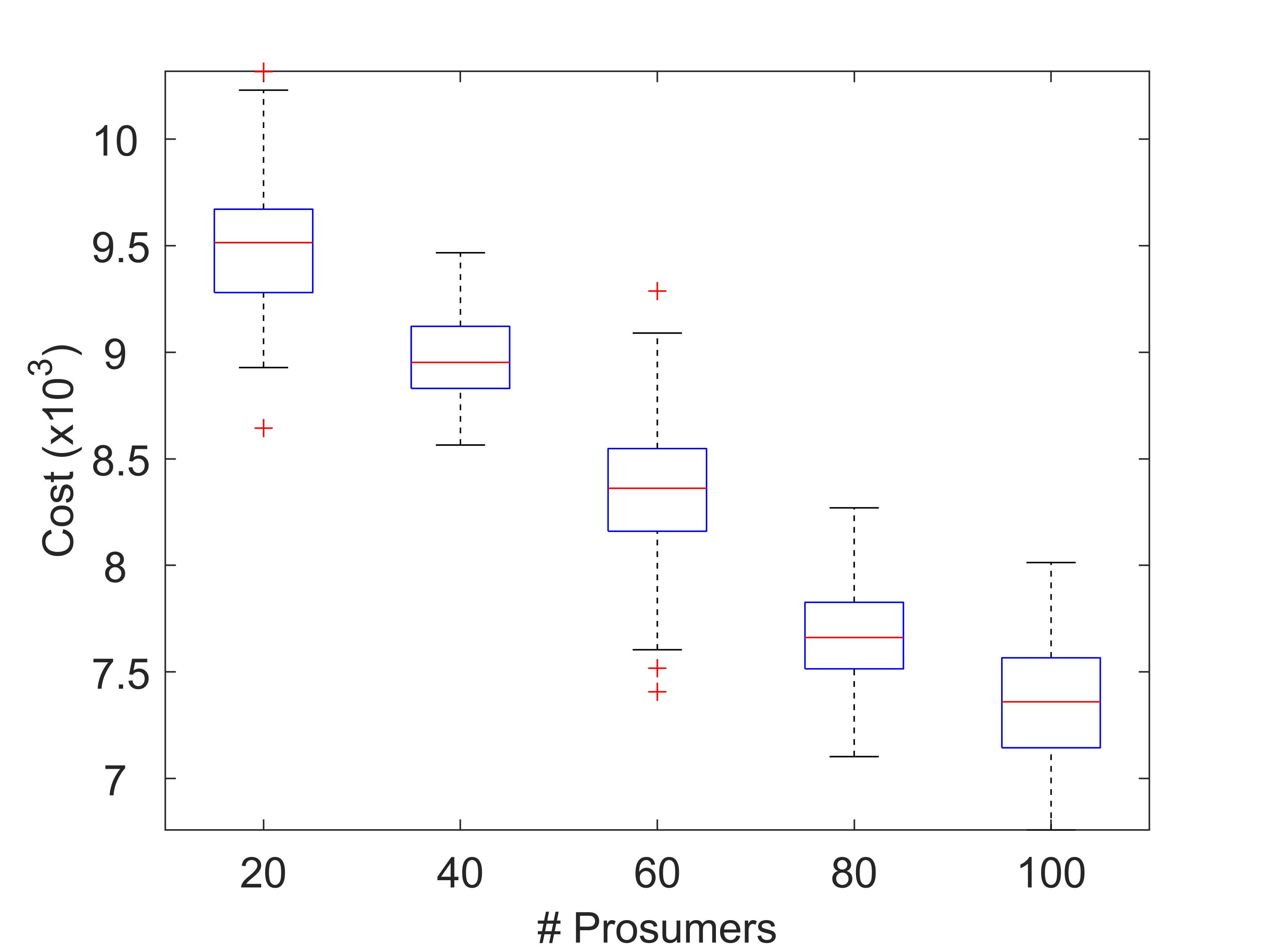}}
		\caption{Impact of model parameters on total cost: (a) Impact of Sample Size, (b) Impact of Prosumer Participation.}
		\label{fig11} 
	\end{figure}

	The proposed algorithm's performance was evaluated by varying the number of prosumers and sample size. As shown in Figure \ref{fig11}a, the total cost decreases with increasing numbers of prosumers. The median total cost consistently drops, and the spread of costs reduces, indicating improved efficiency and stability.
	Figure \ref{fig11}b demonstrates that increasing the sample size leads to a significant reduction in the median total cost. This highlights the importance of larger datasets for accurate uncertainty representation and improved algorithm performance.
	
	\subsubsection{Algorithmic Analysis}
\begin{figure}[!tb]
	\centering
	\subfloat[\label{1a}]{%
		\includegraphics[width=0.5\linewidth]{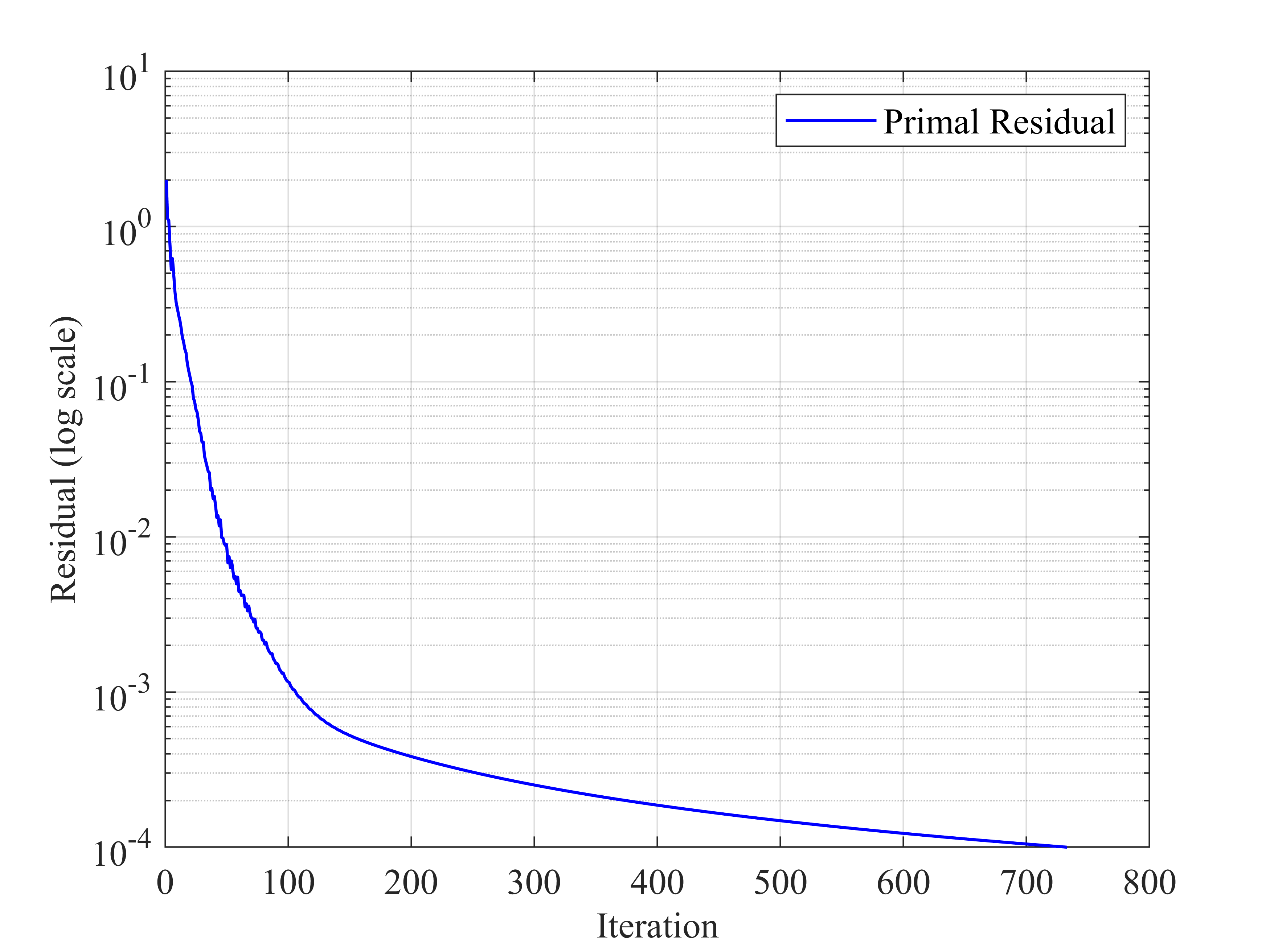}}
	\hfill
	\subfloat[\label{1b}]{%
		\includegraphics[width=0.5\linewidth]{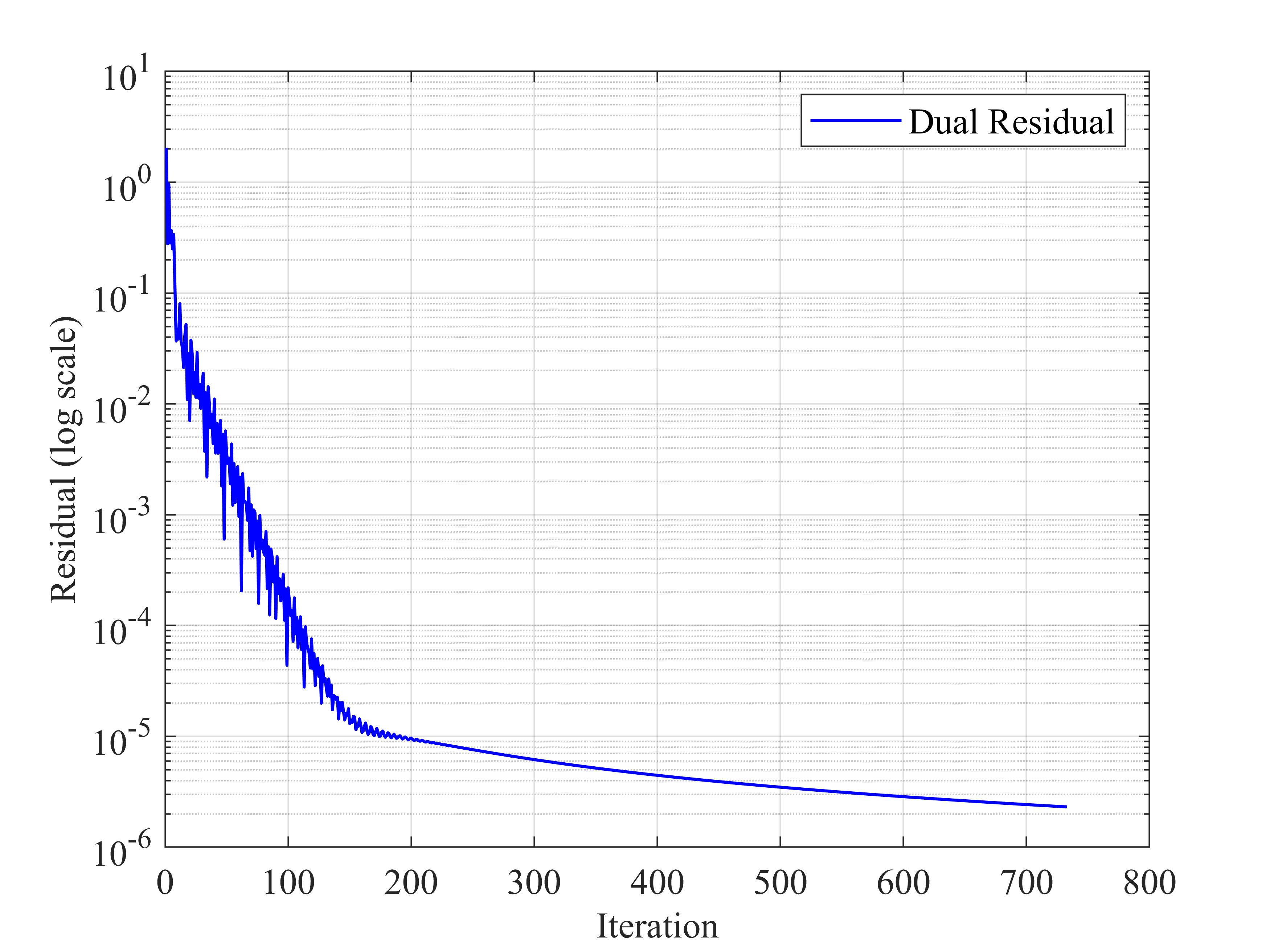}}
	\caption{Residual convergence analysis: (a) Primal residual, (b) Dual residual.}
	\label{fig12} 
\end{figure}

	The convergence behavior of the distributed energy trading algorithm is illustrated in Figure \ref{fig12}. The figure shows the residual over iterations, demonstrating the algorithm's efficiency in achieving optimization objectives.
	Our distributed approach enhances scalability and privacy by limiting information sharing between neighbors. This avoids the communication and computational overhead of centralized solutions, making it suitable for large-scale prosumer networks.

	\section{Conclusion}\label{sec5}
	This paper introduced a novel contextual stochastic optimization framework for coordinated energy trading and flexibility provision, aiming to minimize social costs in distribution networks with a high penetration of renewable energy sources like solar PV. This improvement is partly due to collaborative contextual data exchange between neighboring prosumers. Using a distributed contextual stochastic optimization framework coupled with a distributed optimization algorithm, prosumers can jointly address energy trading and flexibility provision cooperatively, while maintaining privacy by sharing minimal information between neighbors in a computationally efficient and scalable manner.
	Simulation results demonstrated the economic and technical advantages of the proposed approach. By enabling P2P energy trading and coordinating flexibility resources like BESSs and SLs, the model effectively utilized covariates to mitigate output uncertainties from renewables and loads cost-effectively. This framework shows significant potential for facilitating the wide-scale integration of distributed energy resources and the distributed optimization of future smart grids.
	Future work will focus on further validating the model using real system data and expanding the flexibility mechanisms considered for grid services. Additionally, extensions to account for other sources of uncertainties and network constraints are also of interest to make the approach more practical and applicable to distribution system challenges. Overall, the results illustrate the potential of contextual stochastic optimization and distributed algorithms in transitioning to renewable-rich, decentralized energy systems.

	\appendix[KKT conditions for problem \eqref{CKNN} and \eqref{auxUpdate}]
	
	1) \textit{Lower-level problem}: The Lagrangian for the lower level problem is constructed as follows:
	\begin{align*}
		\mathcal{L}_n &= \frac{\gamma}{k} \sum_{j:X_n^j \in \mathcal{N}_k (\hat{x}_{n}^i)} C_n^e(p_n^{mt}, p_{nm}) + C_n^b(p_n^b) + C_n^s(p_n^s) \\
		&+ \frac{1-\gamma}{\lvert \mathcal{N} \rvert_n} \sum_{j:X_n^j \in \mathcal{N}_k (\hat{x}_{mn}^i)} C_n^e(p_n^{mt}, p_{nm}) + C_n^b(p_n^b) + C_n^s(p_n^s) \\
		&+ \mu_p^{\top} (p_n^{s,i} + p_n^{b,i} - p_n^{mt,i} - q_n^{mt,i} - p_n^{g,i} - p_n^{i} ) \\
		&+ \sum_{m\in\mathcal{N}_n} \lambda_{nm} (p_{nm}^{i} + p_{mn}^{i}) \\
		&+ \sum_{h\in\mathcal{H}} \lambda_{b,h} (e_{n,h+1}^i - e_{n,h}^i - \eta_n p_{n,h}^{b,i}) \\
		&+ \sum_{h\in\mathcal{H}} \lambda_{h,t} (h_{n,h+1}^i - h_{n,h}^i - (\hat{p}_{n,h}^{s,i} - \hat{p}_{n,h}^{s,r})) \\
		&+ \underline{\mu}_p^{\top} (\underline{p}_n^{mt} - p_n^{mt,i}) + \overline{\mu}_p^{\top} (p_n^{mt,i} - \overline{p}_n^{mt}) \\
		&+ \underline{\mu}_q^{\top} (\underline{q}_n^{mt} - q_n^{mt,i}) + \overline{\mu}_q^{\top} (q_n^{mt,i} - \overline{q}_n^{mt}) \\
		&+ \underline{\mu}_{pe}^{\top} (\underline{p}_n^{e} - p_n^{e,i}) + \overline{\mu}_{pe}^{\top} (p_n^{e,i} - \overline{p}_n^{e}) \\
		&+ \underline{\mu}_e^{\top} (\underline{e}_n^i-e_n^i) + \overline{\mu}_e^{\top} (e_n^i - \overline{e}_n^i) \\
		&+ \underline{\mu}_b^{\top} (\underline{p}_n^b - p_n^{b,i}) + \overline{\mu}_b^{\top} (p_n^{b,i} - \overline{p}_n^b)\\
		&+ \underline{\mu}_s^{\top} (\underline{p}_n^s - p_n^{s,i}) + \overline{\mu}_s^{\top} (p_n^{s,i} - \overline{p}_n^s) 
	\end{align*}
	By applying the necessary conditions for optimality, the following stationarity conditions are obtained:
	\begin{align*}
		&\frac{\gamma}{k} \sum_{j:X_n^j \in \mathcal{N}_k (\hat{x}_{n}^i)} c_p^{mt}(x_n)
		+ \sum_{j:X_n^j \in \mathcal{N}_k (\hat{x}_{mn})} \frac{1-\gamma}{\lvert \mathcal{N} \rvert_n} c_p^{mt}(x_{nm}) \\
		&- \mu_p^i -\underline{\mu}_p^i + \overline{\mu}_p^i = 0 \\
		&\frac{\gamma}{k} \sum_{j:X_n^j \in \mathcal{N}_k (\hat{x}_{n})} c_q^{mt}(x_n)
		+ \sum_{j:X_n^j \in \mathcal{N}_k (\hat{x}_{mn})} \frac{1-\gamma}{\lvert \mathcal{N} \rvert_n} c_q^{mt}(x_{nm})\\ 
		&- \mu_q^i -\underline{\mu}_q^i + \overline{\mu}_q^i = 0 \\
		&\frac{\gamma}{k} \sum_{j:X_n^j \in \mathcal{N}_k (\hat{x}_{n})} c_{nm}^{e}(x_n)
		+ \sum_{j:X_n^j \in \mathcal{N}_k (\hat{x}_{mn})} \frac{1-\gamma}{\lvert \mathcal{N} \rvert_n} c_{nm}^{e}(x_{nm}) \\ 
		&- \mu_p^i -\underline{\mu}_{pe}^i + \overline{\mu}_{pe}^i + \lambda_{nm} = 0 \\
		&2 \alpha_n^b p_n^b - \underline{\mu_b} + \overline{\mu_b} = 0 \\
		&2 \alpha_n^s p_n^s - \underline{\mu_s} + \overline{\mu_s} = 0 
	\end{align*}
	The primal feasibility conditions are derived from the original constraints of the problem:		
	\begin{align*}		
		&p_n^{l,i} + p_n^{s,i} + p_n^{b,i} - p_n^{mt,i} - q_n^{mt,i} - p_n^{g,i} - P_n^{e,i} \leq  0 \\
		&p_{nm}^{e,i} + p_{mn}^{e,i} = 0 \\
		&e_{n,h+1}^i - e_{n,h}^i - \eta_n p_{n,h}^{b,i} = 0, \quad\quad\quad \forall h\in\mathcal{H} \\
		&s_{n,h+1}^i - s_{n,h}^i - (\hat{p}_{n,h}^{s,i} - \hat{p}_{n,h}^{s,r}) = 0, \,\, \forall h\in\mathcal{H} \\
		&\underline{e}_n^i - e_n^i \leq 0 \\
		&e_n^i - \overline{e}_n^i \leq 0 \\
		&\underline{p}_n^b - p_n^{b,i} \leq 0 \\
		&p_n^{b,i} - \overline{p}_n^b \leq 0 \\
		&\underline{p}_n^s - p_n^{s,i} \leq 0 \\
		&p_n^{s,i} - \overline{p}_n^s \leq 0
	\end{align*}
	The dual feasibility conditions require that the Lagrange multipliers associated with inequality constraints be non-negative:	
	\begin{align*}
		\underline{\mu_p}, \overline{\mu_p}, \underline{\mu_q}, \overline{\mu_q}, \underline{\mu_e}, \overline{\mu_e}, \underline{\mu_b}, \overline{\mu_b}, \underline{\mu_s}, \overline{\mu_s} \geq 0
	\end{align*}
	Finally, the complementary slackness conditions establish a relationship between the Lagrange multipliers and the inequality constraints:	
	\begin{align*}
		&\underline{\mu}_p^{\top} (\underline{p}_n^{mt} - p_n^{mt,i}) = 0 &  &&  && && &&\\
		&\overline{\mu}_p^{\top} (p_n^{mt,i} - \overline{p}_n^{mt}) = 0 \\
		&\underline{\mu}_q^{\top} (\underline{q}_n^{mt} - q_n^{mt,i}) = 0 \\
		&\overline{\mu}_q^{\top} (q_n^{mt,i} - \overline{q}_n^{mt}) = 0 \\
		&\underline{\mu}_{pe}^{\top} (\underline{p}_n^{e} - p_n^{e,i}) = 0 \\
		&\overline{\mu}_{pe}^{\top} (p_n^{e,i} - \overline{p}_n^{e}) = 0 \\
		&\underline{\mu}_e^{\top} (\underline{e}_n^i - e_n^i) = 0 \\
		&\overline{\mu}_e^{\top} (e_n^i - \overline{e}_n^i) = 0 \\
		&\underline{\mu}_b^{\top} (\underline{p}_n^b - p_n^{b,i}) = 0 \\
		&\overline{\mu}_b^{\top} (p_n^{b,i} - \overline{p}_n^b) = 0 \\
		&\underline{\mu}_s^{\top} (\underline{p}_n^s - p_n^{s,i}) = 0 \\
		&\overline{\mu}_s^{\top} (p_n^{s,i} - \overline{p}_n^s) = 0
	\end{align*}
	By employing the bigM method, the nonlinear complementary slackness conditions can be transformed into the following constraints:
	\begin{align*}
		-\underline{\mu}_e &\leq M z_5, & e_n + \overline{e}_n &\leq M (1 - z_5), & -\underline{\mu}_e + e_n + \overline{e}_n &\leq M \\
		\overline{\mu}_e &\leq M z_6, & e_n - \overline{e}_n &\leq M (1 - z_6), & \overline{\mu}_e + e_n - \overline{e}_n &\leq M \\
		-\underline{\mu}_b &\leq M z_3, & p^b_n + \overline{p}^b_n &\leq M (1 - z_3), & -\underline{\mu}_b + p^b_n + \overline{p}^b_n &\leq M \\
		\overline{\mu}_b &\leq M z_4, & p^b_n - \overline{p}^b_n &\leq M (1 - z_4), & \overline{\mu}_b + p^b_n - \overline{p}^b_n &\leq M \\
		-\underline{\mu}_s &\leq M z_1, & s_n + \overline{s}_n &\leq M (1 - z_1), & -\underline{\mu}_s + s_n &\leq M \\
		\overline{\mu}_s &\leq M z_2, & s_n - \overline{s}_n &\leq M (1 - z_2), & \overline{\mu}_s + s_n - \overline{s}_n &\leq M
	\end{align*}

	2) \textit{Auxiliary variable update: }
	To derive the KKT conditions for this problem, the Lagrangian function is constructed as follows: 
	\begin{align*}
		\mathcal{L} = & -(\lambda_{nm,t})^{\top} \hat{e}_{nm,t} - (\lambda_{mn,t})^{\top} \hat{e}_{mn,t} \\
		& + \frac{\rho}{2} \lVert \hat{e}_{nm,t} - e_{nm,t} \rVert_2^2 + \frac{\rho}{2} \lVert \hat{e}_{mn,t} - e_{mn,t} \rVert_2^2 \\
		& + \nu^{\top} (\hat{e}_{nm,t} + \hat{e}_{mn,t}) \\
		& \nabla_{} \mathcal{L} = -\lambda_{nm,t} + \rho (\hat{e}_{nm,t}-e_{nm,t}) + \nu = 0 \\
	\end{align*}
	By applying the necessary conditions for optimality, the following equations are obtained:
	\begin{align*}
		-&\lambda_{nm,t} + \rho (\hat{e}_{nm,t}-e_{nm,t}) + \nu = 0 \\
		-&\lambda_{mn,t} + \rho (\hat{e}_{mn,t}-e_{mn,t}) + \nu = 0 \\
		&\hat{e}_{nm,t}+\hat{e}_{mn,t} = 0
	\end{align*}
	by substituting, the auxiliary variables update rules are derived as 
	\begin{align*}
		& \hat{e}_{nm,t}^{k+1} = \frac{e_{nm,t}^{k+1}-e_{mn,t}^{k+1}}{2} + \frac{\lambda_{nm,t}^{k}-\lambda_{mn,t}^{k}}{2\rho}
	\end{align*}

	\ifCLASSOPTIONcaptionsoff
	\newpage
	\fi

	\bibliographystyle{IEEEtran}
	\bibliography{mybib}

\end{document}